\documentclass[12pt,twoside]{amsart}
\usepackage{amssymb}
\usepackage{amscd}
\usepackage[abbrev,alphabetic]{amsrefs}
\usepackage{hyperref}
\usepackage{comment}
\usepackage{array,multirow,tabularx,longtable}
\usepackage{tikz}
\usetikzlibrary{cd}
\usepackage{here}
\usepackage{multirow}
\usepackage[margin=1.25in]{geometry}

\usepackage{fancybox}
\usepackage{ascmac}

\makeatletter
\@addtoreset{step}{thm}
\makeatother

\makeatletter
\@addtoreset{step}{lem}
\makeatother

\makeatletter
\@addtoreset{step}{prop}
\makeatother

\title[Unboundedness of regular del Pezzo surfaces over imperfect fields]
{Unboundedness of regular del Pezzo surfaces over imperfect fields} 
\author{Hiromu Tanaka} 
\address{Department of Mathematics, 
Graduate School of Science, 
Kyoto University, 
Kyoto 606-8502, JAPAN} 
\email{tanaka.hiromu.7z@kyoto-u.ac.jp}

\newcommand{\Sym}[0]{\operatorname{Sym}}

\newcommand{\NE}[0]{\operatorname{NE}}
\newcommand{\red}[0]{\operatorname{red}}

\newcommand{\Ker}[0]{\operatorname{Ker}}
\newcommand{\Image}[0]{\operatorname{Im}}

\newcommand{\Proj}[0]{\operatorname{Proj}}
\newcommand{\Spec}[0]{\operatorname{Spec}}

\newcommand{\Pic}[0]{\operatorname{Pic}}

\newcommand{\Ex}[0]{\operatorname{Ex}}

\newcommand{\sep}[0]{\operatorname{sep}}

\newtheorem{thm}{Theorem}[section]
\newtheorem{lem}[thm]{Lemma}
\newtheorem{cor}[thm]{Corollary}
\newtheorem{prop}[thm]{Proposition}

\newtheorem*{claim*}{Claim}         
\newtheorem{step}{Step}

\theoremstyle{definition}
\newtheorem{ex}[thm]{Example}

\newtheorem{dfn}[thm]{Definition}

\newtheorem{rem}[thm]{Remark}

\makeatletter
  
  \@addtoreset{equation}{thm}
  \makeatother

\newcommand{\cred}{\color{red}}

\newcommand{\bF}{\mathbb{F}}

\newcommand{\cC}{\mathcal{C}}

\newcommand{\cL}{\mathcal{L}}
\newcommand{\cM}{\mathcal{M}}

\newcommand{\ol}{\overline}

\newcommand{\MO}{\mathcal{O}}

\newcommand{\R}{\mathbb{R}}
\newcommand{\Q}{\mathbb{Q}}
\newcommand{\Z}{\mathbb{Z}}
\newcommand{\F}{\mathbb{F}}

\renewcommand{\P}{\mathbb{P}}

\renewcommand{\P}{\mathbb{P}}
\begin{document}

\maketitle

\begin{abstract}
Over an F-finite field, regular del Pezzo surfaces are bounded. Over a non-F-finite field of characteristic two or three, regular del Pezzo surfaces are unbounded. 
\end{abstract}

\tableofcontents

\section{Introduction}

Fano varieties form one of the central classes of algebraic varieties. 
Their classification and boundedness are among the most fundamental problems in algebraic geometry. 
For example, smooth Fano varieties are classified in dimension at most three \cite{MM81}, \cite{MM83}, \cite{FanoIV}, 
and smooth Fano varieties of fixed dimension form bounded families in characteristic zero \cite{Kol96} (cf.\ \cite{Bir21}).

In this article, we study regular del Pezzo surfaces over imperfect fields of positive characteristic. 
Such surfaces arise naturally as generic fibres of Mori fibre spaces. 
Indeed, let 
\[
\pi:V\to B
\]
be a Mori fibre space over an algebraically closed field of positive characteristic and assume that $V$ has at worst terminal singularities. 
Its generic fibre
\[
X:=V\times_B\Spec K(B)
\]
is defined over the generally imperfect field $K(B)$. 
If $\dim X=2$, then $X$ is a regular del Pezzo surface over $K(B)$. 
Thus, boundedness of regular del Pezzo surfaces over imperfect fields is a natural problem in the study of Mori fibre spaces in positive characteristic.

In \cite{Tan24}, the author established several results related to boundedness of regular del Pezzo surfaces over imperfect fields. 
In particular, 
it is shown  that geometrically integral regular del Pezzo surfaces are bounded. 
The first main result of this article establishes boundedness 
when $[k:k^p]$ is uniformly bounded.

\begin{thm}[Corollary \ref{c F-fin bdd main}]\label{intro bounded}
Fix an integer $e_0>0$. 
Then regular del Pezzo surfaces $X$ over fields $k$ of characteristic $p>0$ satisfying $[k:k^p]\leq p^{e_0}$ are bounded.
\end{thm}

It is then natural to ask whether regular del Pezzo surfaces are uniformly bounded without imposing any restriction on the base field. 
The corresponding question in dimension one has an affirmative answer: every regular Fano curve is isomorphic to a plane conic, and hence regular Fano curves over arbitrary fields form a bounded family.
In dimension two, however, such uniform boundedness fails. 

\begin{thm}[Corollary \ref{c unbounded p=2}, Corollary \ref{c unbounded p=3}]\label{intro unbounded}
For $p \in \{2, 3\}$, 
let $k$ be a field of characteristic $p$ satisfying $[k:k^p]=\infty$. 
Then, for every integer $L>0$, there exists a regular del Pezzo surface $X$ over $k$ such that
\[
K_X^2\geq L
\qquad\text{and}\qquad
h^1(X,\MO_X)\geq L.
\]
\end{thm}

\noindent 
Thus the assumption on $[k:k^p]$ in Theorem \ref{intro bounded} is essential. 
Theorems \ref{intro bounded} and \ref{intro unbounded} together give a sharp contrast between fields of bounded imperfection degree and fields of infinite imperfection degree.

Although we have the above unboundedness result, 
it is still natural to seek the classification of regular del Pezzo surfaces. 
In this direction, 
we establish the following result for the case of index $\geq 2$, 
which plays a crucial role in the proof of Theorem \ref{intro bounded}.

\begin{thm}[Theorem \ref{t index 2}, Corollary \ref{c r>1}]\label{intro index>1}
Let $k$ be a field of characteristic $p>0$ and let $X$ be a projective normal $\Q$-factorial surface over $k$ with $H^0(X, \MO_X)=k$. 
Assume that a positive integer $r$ and an ample Cartier divisor $A$ 
satisfy $-K_X \equiv r A$. 
Then the following hold. 
\begin{enumerate}
\item $r \leq 3$. 
\item If $r =3$, then $X \simeq \P^2_k$. 
\item If $r=2$, then $H^1(X, \MO_X) = H^2(X, \MO_X)=0$ and either 
\begin{itemize}
\item 
$X$ is isomorphic to a quadric surface on $\P^3_k$ 
$($and hence $K_X^2 =8$$)$, or 
\item $p=2$, $\rho(X)=1$, and $K_X^2 =16$. 
\end{itemize}
\end{enumerate}
\end{thm}

\noindent 
There actually exist a regular projective surface $X$ and an ample Cartier divisor $A$ satisfying $-K_X \sim 2A$ and $K_X^2 =16$ (Example \ref{e index2 K^2=16}).

\subsection{Proof of Theorem \ref{intro bounded}}

Fix an integer $e_0>0$. 
Given a field $k$ of characteristic $p>0$ with $[k:k^p] \leq p^{e_0}$ and 
a regular del Pezzo surface $X$ over $k$, 
it is known that the following hold \cite{Tan24}: 
\begin{enumerate}
\item[(a)] $|-12K_X|$ is very ample. 
\item[(b)] $K_X^2 \leq \max \{9, 2^{2e_0+1}\}$.  
\end{enumerate}
By $H^{>0}(X, -12mK_X)=0$ for $m>0$ and the Riemann-Roch theorem, 
the Hilbert polynomial $\chi(X, -12mK_X)$ is given as follows: 
\[
h^0(X, -12mK_X) = \chi(X, -12mK_X) 
= 1-h^1(X, \MO_X)+6m(12m+1) K_X^2. 
\]
In particular, $h^0(X, -12K_X)$ is uniformly bounded, and hence 
the ambient space of the embedding $\varphi_{|-12K_X|} : X \hookrightarrow \P_k^{h^0(X, -12K_X)-1}$ has bounded dimension. 
Then the remaining problem is to control the Hilbert polynomials $\chi(X, -12mK_X)$, 
that is, to find a uniform upper bound of $h^1(X, \MO_X)$. 
In particular, we may assume $H^1(X, \MO_X) \neq 0$. 
Then the problem is reduced to showing the following structure result 
(Theorem \ref{intro structure}), 
because it implies 
\[
h^1(X, \MO_X)-1 = -\chi(X, \MO_X) \leq \frac{1}{2} K_X^2 
\overset{{\rm (b)}}{\leq}\frac{1}{2}\max\{ 9, 2^{2e_0+1}\}. 
\]

\begin{thm}[Theorem \ref{t structure}]\label{intro structure}
Let $X$ be a regular del Pezzo surface satisfying $H^1(X, \MO_X) \neq 0$. 
Then exactly one of the following four cases holds: 
\[
{\setlength{\arraycolsep}{3pt}
\begin{array}{c|c|c|c}
p
& \rho(X)
& \displaystyle-\frac{\chi(X,\MO_X)}{K_X^2}
& \begin{array}{c}\text{Torsor description}\end{array}\\
\hline
2 & 1, 2 & 0
& \begin{array}{c}
\exists\,\alpha_{\omega_X}\text{-torsor}\,\,Y\to X,\\
Y\text{:normal quadric}
\end{array}\\
\hline
2 & 1 & 1/4
& \begin{array}{c}
\exists\,\alpha_{\omega_X}\text{-torsor}\,\,Y\to X,\\
Y\text{:normal},\ K_Y^2=16[k_Y:k]
\end{array}\\
\hline
2 & 1 & 1/2
& \begin{array}{c}
\exists\,A\equiv-2K_X,\\
\exists\,\alpha_{-A}\text{-torsor}\,\,\P^2_{k_Y}\to X
\end{array}\\
\hline
3 & 1 & 1/3
& \begin{array}{c}
\exists\,\alpha_{\omega_X}\text{-torsor}\,\,\P^2_{k_Y}\to X
\end{array}\\
\hline
\end{array}}
\]
\end{thm}

We now give an overview of the proof of Theorem \ref{intro structure}. 
For simplicity, assume that $p=2$ and $H^1(X, 2K_X) \neq 0$. 
It is easy to show $H^1(X, 4K_X)=0$. 
These conditions guarantee the existence of a non-trivial $\alpha_{\omega_X^{\otimes 2}}$-torsor $f: Y \to X$. 
By general results on $\alpha_{\cL}$-torsors, we obtain 
\begin{enumerate}
\item $K_Y \sim f^*(K_X+ (p-1) \cdot 2K_X) = 3f^*K_X$, and 
\item $\chi(Y, \MO_Y) = \sum_{i=0}^{p-1} \chi(X, -2iK_X) = 2\chi(X, \MO_X) + 3K_X^2$. 
\end{enumerate}
By using (1), 
we can prove that $Y \simeq \P^2_{k_Y}$ 
for $k_Y :=H^0(Y, \MO_Y)$. 
Then 
\[
9 [k_Y:k] = K_Y^2 \overset{{\rm (1)}}{=} (3f^*K_X)^2 = (\deg f)\cdot 9 K_X^2,  
\]
which implies $[k_Y:k] = 2K_X^2$. 
On the other hand, it follows from (2) that 
\[
[k_Y:k]=h^0(Y, \MO_Y) = \chi(Y, \MO_Y) = 2\chi(X, \MO_X) + 3K_X^2.
\]
To summarise, we obtain $2K_X^2 = [k_Y:k] = 2\chi(X, \MO_X) + 3K_X^2$, 
which implies $-\chi(X, \MO_X)/K_X^2 = 1/2$. 
The general case can be handled by a similar case-by-case analysis.
For more details, see Section \ref{s bdd}.

\medskip










\subsection{Proof of Theorem \ref{intro unbounded}}

As both cases are similar, we treat only the case when $p=3$. 
Fix an arbitrary integer $n>0$. 
In what follows, we construct a regular del Pezzo surface $X$ satisfying $K_X^2=3^{n+2}$ and $-\chi(X,\MO_X)=3^{n+1}$.

By $[k^{1/3}:k]=[k:k^3]=\infty$, 
we can find an intermediate field
\[
k
\subset
k'
:=
k(\alpha_0,\alpha_1,\alpha_2,\beta_1,\ldots,\beta_n)
\subset
k^{1/3}
\]
such that $[k':k]=3^{n+3}$. 
Set $a_i:=\alpha_i^3\in k$ and $b_r:=\beta_r^3\in k$
for every $0\leq i\leq2$ and $1\leq r\leq n$. 
Set $Y:=\P^2_{k'}=\Proj k'[x_0,x_1,x_2]$.

We define a $k$-derivation
\[
D:k'[x_0,x_1,x_2]
\longrightarrow
k'[x_0,x_1,x_2]
\]
by the following recursive conditions:
\begin{itemize}
\item $D(x_0)=D(x_1)=D(x_2)=0$.
\item $D(\alpha_i)=x_i$ for every $0\leq i\leq2$.
\item $D(\beta_1)=\alpha_0x_1$.
\item
If $1\leq r<n$ and $D(\beta_r)=F_rx_1$ for $F_r\in k'$, then
\[
D(\beta_{r+1})
=
\left(
-\beta_r^2F_r+F_r^3\alpha_1^2
\right)x_1.
\]
\end{itemize}
It is easy to see that $D$ induces the sheaf-theoretic variant:
\[
D:\MO_Y\to\MO_Y(1).
\]
Then we can directly check that $D^3=0$ and $D$ is nowhere vanishing.
By $D(k[x_0,x_1,x_2])=0$, the quotient $f:Y\to Y^D=X$ by $D$ factors as follows:
\[
Y=\P^2_{k'}\xrightarrow{f}X\xrightarrow{g}\P^2_k,
\]
where the composition $Y=\P^2_{k'}\to\P^2_k$ is the base change.
In particular, $\MO_Y(1)=f^*g^*\MO_{\P^2_k}(1)=f^*A$ for $A:=g^*\MO_{\P^2_k}(1)$, which implies that $f:Y\to X$ is an $\alpha_{A^{\otimes-1}}$-torsor (Proposition \ref{p derivation to torsor}).
As $f$ is faithfully flat, $X$ is a regular projective surface.
A direct computation gives $-K_X\equiv A$, $K_X^2=3^{n+2}$, and $-\chi(X,\MO_X)=3^{n+1}$.
It remains to show $H^0(X,\MO_X)=k$.
For $k_r:=k(\alpha_0,\alpha_1,\alpha_2,\beta_1,\ldots,\beta_r)$, 
this can be checked by considering the coefficient derivations
$P_i:k'\to k'$ defined by
$D(\theta)=\sum_{i=0}^2P_i(\theta)x_i$ and their Lie bracket
$C:=[P_0,P_1]$. 
Indeed, one can prove that
$C^{3^{r-1}}|_{k_r}=\partial_{\beta_r}$ for every $1\leq r\leq n$, 
which implies by descending induction that an arbitrary element
$\theta\in H^0(X,\MO_X)=k'^D$ satisfies
$\partial_{\beta_r}(\theta)=0$ and $\partial_{\alpha_i}(\theta)=0$
for every $r$ and $i$.
For more details, see Theorem \ref{t unbounded p=2} ($p=2$) and
Theorem \ref{t unbounded p=3} ($p=3$).

\begin{rem}
There are four types of regular del Pezzo surfaces $X$ with $H^1(X, \MO_X) \neq 0$ listed in Theorem \ref{intro structure}. 
For each type, we shall construct unbounded families (cf.\ Subsection \ref{ss irregular summary}). 
\end{rem}

\subsection{Related results}

In this subsection, we collect some known results related to our main theorems. 


{\em Boundedness and classification.} 
The author established the boundedness of geometrically integral regular del Pezzo surfaces \cite{Tan24}. 
This result was later extended in two directions. 
First, such surfaces were classified by Bernasconi-Tanaka \cite{BT24}. 
In the other direction, 
it is extended to singular settings; 
the boundedness of $\epsilon$-klt del Pezzo surfaces 
is established for the geometrically integral case and the case of characteristic $p>5$ \cite{BM24}, \cite{BD26}. 
Theorem \ref{intro unbounded} shows that this BAB-type result cannot be extended to the case of characteristic two and three.




\medskip

{\em Del Pezzo surfaces with irregularities.} 
The example appearing in Theorem \ref{intro unbounded} 
is a regular del Pezzo surface satisfying $H^1(X, \MO_X) \neq 0$. 
This kind of phenomenon was first discovered by Schr\"{o}er; 
he constructed a weak del Pezzo surface $X$ satisfying $H^1(X, \MO_X) \neq 0$ \cite{Sch07}. 
Then Maddock found the first regular del Pezzo surface $X$ satisfying $H^1(X, \MO_X) \neq 0$ \cite{Mad16}. 
As far as the author knows, Theorem \ref{intro unbounded} 
is the first example that is a regular del Pezzo surface $X$ of characteristic three with $H^1(X, \MO_X) \neq 0$. 
Theorem \ref{intro structure} can be considered as a weak classification result for  regular del Pezzo surfaces with irregularities.




\medskip
\noindent {\bf Acknowledgements.}
The author was supported by JSPS KAKENHI Grant numbers JP23K22383 and JP23K03028. 
The author thanks Fabio Bernasconi for fruitful discussions. 

\medskip
\noindent {\bf AI disclosure.} 
ChatGPT (OpenAI) was used during the preparation of this article for assistance with mathematical exploration, computations, and language editing.
All mathematical statements and proofs were verified by the author, who takes full responsibility for the content.

The results in Sections \ref{s r>1} and \ref{s bdd} were obtained by the author. 
Concerning Sections \ref{s examples} and \ref{s unbounded}, ChatGPT first found the example with 
$r_X=2$ and $K_X^2=16$ (Example \ref{e index2 K^2=16}). 
This example served as the starting point for a series of discussions between the author and ChatGPT, which eventually led to the construction of unbounded examples (Theorem \ref{intro unbounded}).
\section{Preliminaries}

\subsection{Notation}\label{ss-notation}

In this subsection, we summarise notation we will use in this paper. 

\begin{enumerate}
    \item Throughout this article, $p$ denotes a prime number, 
    $k$ denotes a  field of characteristic $p>0$, and 
    we work over $k$ (i.e., we work within the category of $k$-schemes) unless otherwise specified.
\item We will freely use the notation and terminology in \cite{Har77} 
and \cite{Kol13}. 
\item 
For a scheme $X$, its {\em reduced structure} $X_{\red}$ 
is the reduced closed subscheme of $X$ such that the induced morphism 
$X_{\red} \to X$ is surjective. 
\item For an integral scheme $X$, 
we define the {\em function field} $K(X)$ of $X$ 
as $\MO_{X, \xi}$ for the generic point $\xi$ of $X$. 
\item 
We say that $X$ is a {\em variety} (over $k$) or a $k$-{\em variety} if 
$X$ is an integral separated scheme which is of finite type over $k$. 
We say that $X$ is a {\em curve} over $k$ or a $k$-{\em curve} 
(resp. a {\em surface} over $k$ or a $k$-{\em surface}) 
if $X$ is a $k$-variety of dimension one (resp. two). 
\item 
For a variety $X$ over a field $k$, 
its normalisation is denoted by $X^N$. 
For a projective variety $X$, we set $k_X := H^0(X, \MO_X)$. 
\item For a field $k$, 
$\overline k$ (resp. $k^{{\rm sep}}$) 
denotes the algebraic closure (resp. separable closure) of $k$. 
\item For an $\mathbb{F}_p$-scheme $X$ we denote by $F_X \colon X \to X$ the {\em absolute Frobenius morphism}. 
\item If $k \subset k'$ is a field extension and $X$ is a $k$-scheme, we set 
$X_{k'} := X \times_k k' := X \times_{\Spec\,k} \Spec\,k'$. 
\item 
Let $k$ be a field and let $X$ be a normal variety over $k$. 
We say that $X$ is {\em geometrically canonical} 
if $X \times_k \overline k$ is a normal variety over $\overline k$ 
which is canonical, i.e., has at worst canonical singularities. 
Note that if $X$ is geometrically canonical, then $X$ itself is canonical 
\cite[Proposition 2.3]{BT22}. 
\item By abuse of notation, $\Pic X \equiv \Z L$ stands for $(\Pic X/\equiv) = \Z L$. 
\end{enumerate}



\subsection{$\alpha_{\mathcal L}$-torsors}

In this subsection, we collect some known results on $\alpha_{\mathcal L}$-torsors. 
We refer to \cite[Section 1 and Section 2]{Mad16} for more details. 

Given 
a normal variety $X$ and an invertible sheaf $\mathcal L$ on $X$, 
we have an exact sequence on the fppf topology 
\[
0 \to \alpha_{\mathcal L} \to \mathcal L \xrightarrow{F}
\mathcal L^{\otimes p} \to 0,  
\]
which induces an exact sequence 
\begin{equation}\label{e torsor ex seq}
H^0(X, \mathcal L) \xrightarrow{F} H^0(X, \mathcal L^{\otimes p}) \to  H^1_{{\rm fppf}}(X, \alpha_{\mathcal L}) \to H^1(X, \mathcal L) \xrightarrow{F} H^1(X, \mathcal L^{\otimes p}). 
\end{equation}
Here, for every $m \in \Z$, 
$H^1_{{\rm fppf}}(X, \mathcal L^{\otimes m}) \simeq H^1(X, \mathcal L^{\otimes m})$, where the right hand side denotes the usual cohomology with respect to the Zariski topology on $X$.

Note that $F: \mathcal L \to \mathcal L^{\otimes p}$ is a homomorphism of group schemes over $X$ and $\alpha_{\mathcal L}$ is its kernel. 
It is well known that an element $\xi \in H^1_{{\rm fppf}}(X, \alpha_{\mathcal L})$ is bijectively corresponding to 
an $\alpha_{\mathcal L}$-torsor $Y$ over $X$: 
\[
 H^1_{{\rm fppf}}(X, \alpha_{\mathcal L}) 
\overset{{\rm 1:1}}{\leftrightarrow} 
\{ \alpha_{\mathcal L}\text{-torsors over } X\}/\simeq. 
\]
When $Y$ is an $\alpha_{\mathcal L}$-torsor, 
also the induced morphism $f :Y \to X$ is called an {\em $\alpha_{\mathcal L}$-torsor}. 
We say that $\alpha_{\mathcal L}$-torsor is {\em non-trivial} if 
$\xi \neq 0$, i.e., $Y$ is reduced. 
If $L$ is a Cartier divisor on $X$, 
then an $\alpha_L$-torsor means an $\alpha_{\MO_X(L)}$-torsor.  


\begin{rem}\label{r alpha H^1 correspondence}
Let $X$ be a normal variety and let $\mathcal L$ be an invertible sheaf. 
\begin{enumerate}
\item Take nonzero $\zeta \in H^1(X, \mathcal L)$ such that $F(\zeta)=0$ in $H^1(X, \mathcal L^{\otimes p})$. 
By the exact sequence (\ref{e torsor ex seq}), there exists a nonzero element $\xi \in H^1_{{\rm fppf}}(X, \alpha_{\mathcal L})$ mapping to $\zeta$. 
Therefore, such an element $\zeta$ induces an $\alpha_{\mathcal L}$-torsor. 
\item 
If $H^0(X, \mathcal L^{\otimes p})=0$ 
(we will only encounter this situation in this article), then we get 
\[
\Ker(F : H^1(X, \mathcal L) \to H^1(X, \mathcal L^{\otimes p}))
= H^1_{{\rm fppf}}(X, \alpha_{\mathcal L}) 
\overset{{\rm 1:1}}{\leftrightarrow} 
\{ \alpha_{\mathcal L}\text{-torsors over }X\}/\simeq
\]
\item 
By (2), 
if $H^0(X, \mathcal L^{\otimes p})=0$ and there exists an $\alpha_{\mathcal L}$-torsor $Y \to X$ such that $Y$ is reduced, 
then $H^1(X, \mathcal L) \neq 0$. 
\end{enumerate}
\end{rem}


\begin{prop}\label{p torsor}
Let $X$ be a projective normal Gorenstein variety. 
Take an invertible sheaf $\mathcal L$  and a non-trivial $\alpha_{\mathcal L}$-torsor $f: Y \to X$. 
Then the following hold. 
\begin{enumerate}
\item $Y$ is a projective Gorenstein variety over $k$. 
\item $\omega_Y \simeq f^*(\omega_X \otimes \mathcal L^{p-1})$. 
\item $\chi(Y, \MO_Y) = \sum_{i=0}^{p-1} \chi(X, \mathcal L^{-i})$, 
where $\chi(Y, \MO_Y) := \sum_{j \geq 0} (-1)^j\dim_k H^j(Y, \MO_Y)$ 
and $\chi(X, \mathcal L^{-i}) := \sum_{j \geq 0} (-1)^j\dim_k H^j(X, \mathcal L^{-i})$. 
\end{enumerate}
\end{prop}

\begin{proof}
See \cite[Proposition 1.1.2]{Mad16}. 
\end{proof}

\begin{prop}\label{p torsor2}
Assume that $p=2$. 
Let $f: Y \to X$ be a finite surjective morphism  
from a CM projective variety $Y$ to a regular projective variety $X$ 
such that 
the induced field extension $K(Y)/K(X)$ is inseparable and of degree two. 
Then $f: Y \to X$ is a non-trivial $\alpha_{\mathcal L}$-torsor for some invertible sheaf $\mathcal L$ on $X$. 
\end{prop}

\begin{proof}
See \cite[Proposition 2.2.1]{Mad16}. 
\end{proof}

\begin{rem}
Although we cite \cite[Section 1 and Section 2]{Mad16} as a reference on $\alpha_{\mathcal L}$-torsors, 
we need to be careful with it, because 
\cite[Lemma 1.1.3]{Mad16} is incorrect 
(indeed, $\varphi : Y \to X$ in \cite[Proposition 9.8]{Tan21} is a finite surjective morphism of degree $p$ 
between regular projective varieties, whilst $[H^0(Y, \MO_Y) : H^0(X, \MO_X) ] =p^2$). 
This error affects  
\[
\text{Theorem 1.2.3,\quad Corollary 1.2.6,\quad  Corollary 1.2.8,\quad 
Corollary 2.2.2}
\]
of \cite{Mad16}. 
\end{rem}




\subsection{Del Pezzo surfaces}\label{ss reg dP prelim}

A {\em del Pezzo surface} $X$  
is a projective normal surface such that $H^0(X, \MO_X)=k$ and 
$-K_X$ is an ample $\Q$-Cartier divisor. 
In particular, a regular del Pezzo surface $X$ is a regular projective surface such that 
$H^0(X, \MO_X)=k$ and $-K_X$ is ample.

\begin{prop}\label{p dP basic1}
Let $X$ be a regular del Pezzo surface. 
Then the following hold. 
\begin{enumerate}
\item $X$ is geometrically normal if and only if $X$ is geometrically canonical. 
\item If $p \geq 5$, then $X$ is geometrically canonical. 
\end{enumerate}
\end{prop}

\begin{proof}
The assertion (1) (resp. (2)) 
follows from \cite[Theorem 3.3]{BT22} 
(resp. (1) and \cite[Theorem 1.5]{PW22}). 
\end{proof}

\begin{prop}\label{p dP basic2}
Let $X$ be a regular del Pezzo surface. 
Then the following hold. 
\begin{enumerate}
\item If $X$ is geometrically normal, then $H^i(X, K_X+A)=0$ for every 
$i>0$ and ample Cartier divisor $A$ $($in particular, $H^1(X, \MO_X)=0)$.
\item If $H^1(X, \MO_X) = 0$, then $\Pic\,X \simeq \Z^{\oplus \rho(X)}$ 
and every numerically trivial Cartier divisor $N$ on $X$ satisfies $N \sim 0$. 
\end{enumerate}
\end{prop}

\begin{proof}
Let us show (1). 
The following holds for $j <2$ and the pullback $A_{\overline k}$ of $A$ on $X_{\overline k}$: 
\[
H^j(X, -A) \otimes_k \overline k 
\simeq 
H^j(X_{\overline k}, -A_{\overline k}) 
\overset{{\rm (\alpha)}}{\simeq} H^j(Y, -\mu^*A_{\overline k}) 
\overset{{\rm (\beta)}}{=} 0, 
\]
where $\mu : Y \to X_{\overline k}$ is the minimal resolution, 
$(\alpha)$ follows from the fact that $X_{\overline k}$ has at worst rational singularities, and 
($\beta$) holds by \cite{Muk13}. 
Thus (1) holds.

Let us show (2). 
Take a numerically trivial Cartier divisor $N$ on $X$. 
It is enough to prove that $H^0(X, N) \neq 0$. 
By the Riemann-Roch theorem, we get 
\[
h^0(X, N) - h^1(X, N) = \chi(X, N) = \chi(X, \MO_X) = h^0(X, \MO_X) =1. 
\]
Thus (2) holds. 
\end{proof}

\section{Case $r_X \geq 2$}\label{s r>1}

\begin{dfn}
Let $X$ be a regular del Pezzo surface. 
The {\em index} $r_X$ of $X$ is defined as the largest positive integer $s$ such that there exists an ample Cartier divisor $A$ on $X$ satisfying $-K_X \equiv sA$. 
\end{dfn}

In this section, we treat the case $r_X \geq 2$. 
Specifically, we prove that if $X$ is a regular del Pezzo surface 
with $r_X \geq 2$, then the following hold (Corollary \ref{c r>1}). 
\begin{enumerate}
\item $r_X \in \{2, 3\}$. 
\item If $r_X =3$, then $X \simeq \P^2_k$. 
\item If $r_X =2$, then one of the following holds: 
\begin{itemize}
\item $K_X^2 = 8$ and 
$X$ is isomorphic to a quadric surface on $\P^3_k$. 
\item $K_X^2 = 16$ and $H^1(X, -A) \neq 0$ for some ample Cartier divisor $A$. 
\end{itemize}
\end{enumerate}
For later usage, we will establish these results for the case when $X$ is singular (more precisely, 
the case when $X$ is a projective normal $\Q$-factorial surface). 
Such a generalisation is necessary, because we shall apply this result to 
some $\alpha_{\mathcal L}$-torsor of a given regular del Pezzo surface, which is not necessarily regular. 

\subsection{Case $r_X =3$}

\begin{lem}\label{l higher index}
Let $X$ be a projective normal $\Q$-factorial surface with $H^0(X, \MO_X)=k$. 
Let $\mu : Y \to X$ be the minimal resolution of $X$. 
Assume that there exist a nef $\Q$-divisor $N$ on $X$ 
and ample Cartier divisors $A$ and $A'$ on $X$ 
such that $-K_X \equiv A +A' +N$. 
Then the following hold. 
\begin{enumerate}
\item $\rho(X) \leq 2$. 
\item If $\rho(X) =2$, 
then $X$ is regular and 
the contraction of each extremal ray 
is a Mori fibre space to a curve. 
\item 
If $\rho(X) =1$ and $X$ is not regular, then $X$ is canonical, $\rho(Y)=2$, 
and $Y$ has a $K_Y$-Mori fibre space to a curve. 
\item Every effective divisor on $X$ is nef. 
\end{enumerate}
In particular, $X$ is canonical. 
\end{lem}

\begin{proof}
If $\rho(Y) =1$ (i.e., $X$ is regular and $\rho(X)=1$), then there is nothing to show. 
In what follows, we assume $\rho(Y) \geq 2$. 
We have $K_Y+\Delta  \sim_{\Q} \mu^*K_X$ for some effective $\mu$-exceptional $\Q$-divisor $\Delta$, which implies 
\[
K_Y + \Delta +  \mu^*A + \mu^*A' +\mu^*N \equiv 0. 
\]

We now prove that 
\begin{enumerate}
\item[(1)'] 
If $R = \R_{\geq 0}[C]$ is a $K_Y$-negative extremal ray of $\overline{\NE}(Y)$ 
and $C$ is a prime divisor on $Y$, then 
$C^2  =0$, $\rho(Y) =2$, and 
the contraction $\pi : Y \to B$ of $R$ 
is a $K_Y$-negative Mori fibre space to a curve $B$. 
\end{enumerate}
Take a $K_Y$-negative extremal ray  $R = \R_{\geq 0}[C]$ of $\overline{\NE}(Y)$, 
where $C$ is a prime divisor on $Y$. 
By $\rho(Y) \geq 2$, we get $C^2 \leq 0$. 
Suppose $C^2 <0$. We then get the following contradiction: 
\[
0 =(K_Y + \Delta + \mu^*A + \mu^*A' +\mu^*N) \cdot C 
\geq h^0(C, \MO_C) (-1 + 0 + 1 + 1 +0)= h^0(C, \MO_C), 
\]
where $h^0(C, \MO_C) := \dim_k H^0(C, \MO_C)$. 
Hence $C^2 = 0$, which completes the proof of (1)'. 
In particular, we get 
$\rho(X) \leq \rho(Y) =2$, and hence (1) holds. 


Let us show (2). 
Assume $\rho(X)=2$. 
This, together with $\rho(X) \leq \rho(Y) =2$, implies $\rho(X) = \rho(Y)$, i.e., $X = Y$. 
Then the assertion (2) follows from (1)'. 
Note that (4) follows from (1) and (2).

It is enough to show (3). 
Assume $\rho(X)=1$ and $X$ is not regular. 
By (1)', $\rho(Y)=2$ and 
we have a $K_Y$-negative Mori fibre space $\pi : Y \to B$ 
to a regular projective curve $B$. 
Hence it suffices to prove that $X$ is canonical. 
Suppose that $X$ is not canonical, i.e., $\Delta \neq 0$. 
For the $\mu$-exceptional prime divisor $E: = \Ex(\mu)$, 
we have $\Delta = b E$ for some rational number $b>0$. 
Since $E$ is a nonzero effective $\Q$-divisor 
(and hence $-E$ is not nef), 
we get $\deg_{k_F}(\Delta|_F) = \deg_{k_F}(b E|_F) >0$ for the generic fibre $F$ of $\pi : Y \to B$ (note that $\NE(Y) = \R_{\geq 0}[E]+\R_{\geq 0}[F]$). 
We then get the following contradiction: 
\[
0 = \deg_{k_F} ((K_Y+\Delta -\mu^*K_X)|_F) 
\]
\[
= 
\deg_{k_F} (K_F + (\Delta +\mu^*A + \mu^*A' +\mu^*N)|_F ) > -2 + 0 +1+1+0 = 0.
\]
Thus (3) holds. 
\qedhere

\end{proof}

\begin{prop}\label{p index 3}
Let $X$ be a projective normal $\Q$-factorial surface with $H^0(X, \MO_X)=k$. 
Assume that there exist 
a nef $\Q$-divisor $N$ and 
ample Cartier divisors $A$ 
and $A'$
such that $-K_X \equiv A + A' +N$ and $N \not\equiv 0$. 
Then $X \simeq \P^2_k$ and  $A \sim A' \sim \MO_{\P^2}(1)$. 
\end{prop}


\begin{proof}
We start by showing the following claim.

\begin{claim*}
$X$ is regular and $\rho(X)=1$. 
\end{claim*}

\begin{proof}[Proof of Claim]
Let $\mu : Y \to X$ be the minimal resolution. 
Suppose $\rho(Y) \geq 2$. It suffices to derive a contradiction. 
By Lemma \ref{l higher index}, 
$\rho(Y)=2$ and there is a $K_Y$-negative Mori fibre space $\pi: Y \to B$ to a regular projective curve $B$. 

Suppose $\rho(X) =1$. 
Then $N$ is ample and we have $K_Y+\Delta  \sim_{\Q} \mu^*K_X$ for some $\mu$-exceptional effective $\Q$-divisor $\Delta$ on $Y$, which implies 
\[
K_Y + \Delta + \mu^*A  + \mu^*A' + \mu^*N \equiv 0. 
\]
We get the following contradiction 
for a fibre $F$ of $\pi : Y \to B$: 
\[
0 =(K_Y + \Delta + \mu^*A + \mu^*A' + \mu^*N ) \cdot F
> [k_F:k](-2 + 0 + 1+1 + 0)=0. 
\]

Hence we get $\rho(X) >1$. 
By $2 \geq \rho(Y) \geq \rho(X) >1$, we have $\rho(X)=2$. 
By Lemma \ref{l higher index}(2), 
both the extremal rays of $X$ induce Mori fibre space structures to curves, i.e., 
$X$ has two Mori fibre space structures $\pi_i : X \to B_i$, 
where $B_i$ is a regular projective curve for each $i \in \{1, 2\}$. 
It follows from $N \not\equiv 0$ that $N \cdot F_j>0$ for some $j\in \{1, 2\}$, 
where  $F_j$ is a fibre of $\pi_j$. 
This leads to the following contradiction: 
\[
0 = (K_Y + \Delta + \mu^*A+ \mu^*A' + \mu^*N ) \cdot F_j 
>  [k_{F_j}:k](-2 +0 +1 +1 + 0)=0. 
\]
This completes the proof of Claim. 
\end{proof}


If $k$ is an algebraically closed field, 
then the assertion immediately follows from Claim, 
because a smooth del Pezzo surface of Picard number one over an algebraically closed field is isomorphic to $\P^2$. 
In what follows, we reduce the general case to this case. 

Let 
\[
f : Z := (X \times_k \overline k)^N_{\red} \to X
\]
be the induced morphism. 
Then \cite[Theorem 4.2]{Tan18b} enables us to find an effective Weil divisor $C$ on $Z$ 
such that 
\begin{itemize}
\item $K_Z + C \sim f^*K_X$, and 
\item $C =0$ if and only if $X$ is geometrically normal. 
\end{itemize}
In particular, 
\[
K_Z + C + f^*A+f^*A' +f^*N \equiv 0. 
\]
By $\rho(Z) = \rho(X)=1$ \cite[Lemma 2.2 and Proposition 2.4]{Tan18b}, 
$C$ is nef. 
Since Proposition \ref{p index 3} has been already settled 
when $k$ is algebraically closed, 
we get that $Z \simeq \P^2_{\overline k}$ and $C + f^*N \sim \MO_{\P^2}(1)$. 
As $C$ is an effective Weil (so Cartier) divisor and $f^*N$ is a nef $\Q$-Cartier $\Q$-divisor with $f^*N \not\equiv 0$, 
we get $C =0$ and $f^*N \sim \MO_{\P^2}(1)$. 
Therefore, $X$ is geometrically normal, i.e., $X \times_k \overline k \simeq Z \simeq \P^2_{\overline k}$. 
Since $|A_{\overline k}|$ is very ample and 
$h^0(X_{\overline k}, A_{\overline k}) = 3$, 
also $|A|$ is very ample and $h^0(X, A)=3$. 
Hence we get a closed immersion 
$\varphi_{|A|} : X \hookrightarrow \P^2_k$, which is an isomorphism. 
\end{proof}

\subsection{Case $r_X =2$}

\begin{lem}\label{l P^1 bundle}
Let $X$ be a regular projective surface  with $H^0(X, \MO_X)=k$ and     
let $\pi : X \to B$ be a $K_X$-negative Mori fibre space to 
a regular projective curve $B$. 
Assume that there exist $\pi$-ample Cartier divisors $A$ and $A'$
on $X$ satisfying 
$-K_X \equiv_{\pi} A+A'$. 
Then $\pi_*\MO_X(A)$ is a locally free sheaf of rank $2$, 
$X \simeq \P_B(\pi_*\MO_X(A))$, and 
$\pi$ coincides with the canonical projection. 
Moreover, $K_X^2 = 8(1-h^1(B, \MO_B)).$ 
\end{lem}

\begin{proof}
Fix a (not necessarily closed) point $b \in B$. 
Then the scheme-theoretic fibre $X_b$ of $\pi$ over $b$ 
is an integral scheme which is a conic on $\P^2_{k_b}$ \cite[Proposition 2.18]{BT22}. 
By $-K_X \equiv_{\pi} A + A'$  and \cite[Lemma 10.6]{Kol13}, 
we get $X_b \simeq \P^1_{k_b}$ and $A|_{X_b} \sim \MO_{\P^1}(1)$. 
Therefore, we get $\dim_{k(b)} H^0(X_b, A|_{X_b}) = 2$, which implies that $\pi_*\MO_X(A)$ is locally free sheaf of rank $2$. 
It is easy to see that  $\MO_X(A)$ is relatively very ample over $B$. 
Hence we get a closed immersion 
\[
X \hookrightarrow \P_B(\pi_*\MO_X(A)), 
\]
which must be an isomorphism by $\dim X = \dim \P_B(\pi_*\MO_X(A))$. 
The remaining equality  $K_X^2 = 8(1-h^1(B, \MO_B))$ 
follows from the same argument as in \cite[Lemma 4.8]{BT24}. 
\end{proof}

\begin{prop}\label{p index 2 pic2}
Let $X$ be a projective normal $\Q$-factorial surface with $H^0(X, \MO_X)=k$. 
Assume that $\rho(X) \geq 2$ and $-K_X \equiv A+A'$ for some ample Cartier divisors $A$ and $A'$ on $X$. 
Then it holds that 
\[
X \simeq B_1 \times_k B_2,  
\]
where $B_1$ and $B_2$ are smooth conics on $\P^2_k$. 
\end{prop}

\begin{proof}



We first reduce the problem to the case when $k$ is separably closed. 
By Lemma \ref{l higher index}, $X \times_k k^{{\rm sep}}$ is canonical and hence the assumptions hold also for $X \times_k k^{{\rm sep}}$ \cite[Corollary 4.11]{Tan18m}. 
Again by Lemma \ref{l higher index}, both $X$ and $X \times_k k^{{\rm sep}}$  are regular and  $\rho(X) = \rho(X \times_k k^{{\rm sep}})=2$. 
Moreover, we have 
two Mori fibre spaces 
$\pi_1 : X \to B_1$ and $\pi_2: X \to B_2$ to curves $B_1$ and $B_2$. 
Set $\pi := \pi_1 \times \pi_2 : X \to B_1 \times_k B_2$ to be the induced morphism. 
Then $\pi$ is an isomorphism, because so is $\pi \times_k k^{\sep}$ 
(as we are assuming that the assertion holds for the case when $k$ is separably closed).

\medskip

Suppose that $X$ is not geometrically normal. 
We can find a finite purely inseparable extension $l/k$ 
such that $X_l := X \times_k l$ is an integral scheme which is not normal 
by \cite[Lemma 1.3]{Sch10} or \cite[Lemma 3.3]{Tan18b}. 
For the normalisation $Y$ of $X_l$, 
we have the induced morphisms $f: Y  \to X_l \to X$ and 
\[
K_Y +C \equiv f^*K_X \equiv -f^*(A+A'), 
\]
where $C$ is a nonzero effective Weil divisor on $Y$. 
Hence $-K_Y \equiv f^*A +f^*A' +C$. 
Recall that every effective divisor on $X$ is nef, 
and hence the same property holds for $Y$. 
Then $C$ is a nonzero nef $\Q$-Cartier Weil divisor. 
It follows from Proposition \ref{p index 3} that $Y \simeq \P^2$, 
which leads to a contradiction: $1 = \rho(Y) = \rho(X) \geq 2$. 
Therefore, $X$ is geometrically normal.

Recall that we have a Mori fibre space $\pi : X \to B$ 
to a curve $B$ (Lemma \ref{l higher index}). 
By $-K_X \equiv A+A'$, 
$\pi : X \to B$ is a $\P^1$-bundle (Lemma \ref{l P^1 bundle}), and hence $-K_B$ is ample \cite[Corollary 4.9(2)]{Eji19}. 
Since $X$ is geometrically normal, so is $B$. 
Then $B$ is a smooth conic. 
As $k$ is separably closed, we get $B \simeq \P^1_k$. 
Hence $X =\P_{\P^1}(E)$ for some locally free sheaf $E$ of rank $2$. 
We may assume that $E \simeq \MO_{\P^1} \oplus \MO_{\P^1}(n)$ for some $n \in \Z$. 
Since $X$ has two Mori fibre spaces to curves (Lemma \ref{l higher index}), 
we get $n=0$, i.e., $X \simeq \P^1 \times \P^1$. 
\qedhere



\end{proof}

\begin{thm}\label{t index 2}
Let $X$ be a projective normal $\Q$-factorial surface with $H^0(X, \MO_X)=k$. 
Assume that $\rho(X) =1$ and 
there exists an ample Cartier divisor $A$ such that $-K_X \equiv 2A$. 
Then the following hold: 
\begin{enumerate}
\item $X$ is canonical. 
\item $\Pic X \equiv \Z A$. 
\item $H^1(X, M)=0$ if $M$ is a Cartier divisor satisfying $M \not\equiv -A$. 
\item 
One of the following holds: 
\begin{enumerate}
\item $K_X^2 =8$, $H^1(X, \MO_X(-A))=0$, and $X$ is isomorphic to a quadric surface on $\P^3_k$. 
Moreover, $\Pic\,X = \Z A$. 
\item $p=2$, $K_X^2=16$, and $h^1(X, \MO_X(-A)) =1$. 
\end{enumerate}
\end{enumerate}

\end{thm}

\begin{proof}

\setcounter{step}{0}




\begin{step}\label{s1 index 2}
The following hold. 
\begin{enumerate}
\renewcommand{\labelenumi}{(\roman{enumi})}
\item $X$ is canonical. 
\item $\Pic\,X \equiv \Z A$. 
    \item If $m \in \Z$ and $M$ is a Cartier divisor satisfying $M \equiv mA$ and $m \neq -1$, then $H^1(X, M)=0$. 
    In particular, $H^1(X, \MO_X) =0$.
    \item If $H^1(X, \MO_X(-A))=0$, then $|A|$ is very ample.        
\end{enumerate}
\end{step}

\begin{proof}[Proof of Step \ref{s1 index 2}]
The assertion (i) follows from Lemma \ref{l higher index}. 
By 
$\rho(X)=1$, 
we have $\Pic\,X \equiv \Z H$ for some ample Cartier divisor $H$ on $X$. 
In particular, $A \equiv bH$ for some integer $b>0$, 
which implies $-K_X \equiv 2A \equiv 2bH$. 
By Proposition \ref{p index 3}, we get $b=1$, i.e., (ii) holds. 


Let us show (iii). By Serre duality, 
we may assume $m \leq -2$. 
Set $n:= -m \geq 2$, and hence $M \equiv -nA$. 
Suppose that 
$H^1(X, M) \neq 0$. 
By replacing $M$ by $p^eM$ for some integer $e\geq 0$ (if necessary), 
we may assume that 
$H^1(X, M) \neq 0$ and $H^1(X, pM)=0$, 
which is guaranteed by the Serre vanishing theorem. 
Then 
a nonzero element $\eta \in H^1(X, M) \setminus \{0\}$ 
induces an $\alpha_{M}$-torsor $\psi : X' \to X$ 
(Remark \ref{r alpha H^1 correspondence}), 
which satisfies the following properties (Proposition \ref{p torsor}):  
\begin{itemize}
\item $\psi$ is a finite purely inseparable surjective morphism of degree $p$, 
\item $X'$ is a projective Gorenstein surface, and 
\item $K_{X'} \sim \psi^*(K_X + (p-1)M) \equiv -(2 + (p-1)n)\psi^*A$. 
\end{itemize}
Taking the normalisation $X'' \to X'$ of $X'$, 
it holds that $K_{X''} +C \equiv -(2 + (p-1)n)A''$ for some ample Cartier divisor $A''$ and an effective divisor $C$. 
By $\rho(X) = 1$, we get $\rho(X'') =1$ \cite[Proposition 2.4(1)]{Tan18b}, and hence $C$ is nef. 
However, this is absurd by $2 + (p-1)n \geq 2 + 2 >3$ (Proposition \ref{p index 3}). 
Thus (iii) holds. 

Let us show (iv). 
We have 
\[
h^0(X, A) \geq \chi(X, A) = 1 + \frac{1}{2}A \cdot (A -K_X) >0, 
\]
which implies $H^0(X, A) \neq 0$. 
Pick $\Gamma \in |A|$, which is an effective Cartier divisor. 
By the exact sequence $0 \to \MO_X(-A) \to \MO_X \to \MO_{\Gamma} \to 0$, 
we get $H^1(\Gamma, \MO_{\Gamma})=0$ and $H^0(X, \MO_X) \xrightarrow{\simeq} H^0(\Gamma, \MO_{\Gamma})$. 
For an integer $d >0$, we get the following exact sequence by (iii): 
\[
H^0(X, \MO_X(dA)) \to H^0(\Gamma, \MO_X(dA)|_{\Gamma}) \to H^1(X, \MO_X(dA -\Gamma))=0. 
\]
By \cite[Proposition 3.2(3)]{Tan25}, 
$\bigoplus_{d=0}^{\infty} H^0(X, \MO_X(dA))$ is generated by 
$H^0(X, \MO_X(A))$ as a $k$-algebra (cf. \cite[Theorem (2.3)]{Fuj90}). 
In particular, $|A|$ is very ample. 
Thus (iv) holds. 
This completes the proof of Step \ref{s1 index 2}. 
\end{proof}

\begin{step}\label{s2 index 2}
If $H^1(X, -A) =0$, 
then $X$ is isomorphic to a quadric surface on $\P^3_k$. 
\end{step}

\begin{proof}[Proof of Step \ref{s2 index 2}]
By the Riemann-Roch theorem, we get 
\[
0 = \chi(X, -A) = 1 + \frac{1}{2}(-A) \cdot (-A -K_X) = 1 -\frac{1}{2} A^2, 
\]
which implies $A^2 =2$. 
This, together with (iii), implies 
\[
h^0(X, A) = \chi(X, A) = 1 + \frac{1}{2}A \cdot (A -K_X) =1 + \frac{3}{2}A^2 = 4. 
\]
By (iv), $X$ is (isomorphic to) a surface on $\P^3_k$, 
which is a quadric surface by $-K_X \equiv 2A$. 
This completes the proof of Step \ref{s2 index 2}. 
\end{proof}

\begin{step}\label{s3 index 2}
If $H^1(X, -A) \neq 0$, then $A^2 =4$, $K_X^2 =16$, and 
$h^1(X, -A) =1$.  
\end{step}

\begin{proof}[Proof of Step \ref{s3 index 2}]
Assume $H^1(X, -A) \neq 0$. 
By $H^1(X, -pA)=0$ and Remark \ref{r alpha H^1 correspondence}, 
we can find a non-trivial $\alpha_{-A}$-torsor $\psi : X' \to X$. 
Then $X'$ is a projective Gorenstein surface satisfying $K_{X'} \sim \psi^*(K_X - (p-1)A)$ (Proposition \ref{p torsor}). 
For the normalisation $\nu : X'' \to X'$ of $X'$, we get 
\[
K_{X''} +C = \nu^*K_{X'} \sim \nu^*\psi^*(K_X - (p-1)A)
\]
for the conductor divisor $C$ of $\nu : X'' \to X'$. 
Then it holds that  
\[
-K_{X''} \equiv  C + (p+1) \nu^*\psi^*A.  
\]
By $\rho(X'') =\rho(X) =1$, $C$ is nef. Then it follows from Proposition \ref{p index 3} that 
$p=2$, $C=0$, $X' = X'' = \P^2_{k'}$, and 
$-K_{X'} \sim 3\psi^*A$, where $k' := H^0(X', \MO_{X'})$. 
We then get 
\[
9[k':k] = K_{X'}^2 = (3\psi^*A)^2 = 9 \cdot (\deg \psi) \cdot A^2 = 9 \cdot 2 \cdot A^2, 
\]
and hence $[k': k] = 2A^2$. 
By Proposition \ref{p torsor}(3), we get 
\[
[k':k] = \chi(X', \MO_{X'}) = \chi(X, \MO_X) +\chi(X, A) = 1 + (1 + \frac{1}{2}A \cdot (A -K_X)) = 2 + \frac{3}{2} A^2. 
\]
Therefore, we get 
$A^2 =4$ and $[k':k]= 8.$ 

We have $H^0(X, -A) = H^2(X, -A)=0$. 
Hence the Riemann--Roch theorem implies 
\[
-h^1(X,-A)
=
\chi(X,-A)
=
1+\frac{1}{2}(-A)\cdot(-A-K_X)
=
1-\frac{1}{2}A^2
=
-1.
\]
This completes the proof of Step \ref{s3 index 2}. 
\qedhere
\end{proof}
This completes the proof of Theorem \ref{t index 2}. 
\end{proof}

\begin{cor}\label{c r>1}
Let $X$ be a projective normal $\Q$-factorial surface with $H^0(X, \MO_X)=k$. 
Assume that a positive integer $r$ and an ample Cartier divisor $A$ 
satisfy $-K_X \equiv r A$. 
Then the following hold. 
\begin{enumerate}
\item $r \leq 3$. 
\item If $r =3$, then $X \simeq \P^2_k$. 
\item If $r=2$, then $H^1(X, \MO_X) = H^2(X, \MO_X)=0$ and either 
\begin{itemize}
\item 
$X$ is isomorphic to a quadric surface on $\P^3_k$ (and hence $K_X^2 =8$), or 
\item $p=2$, $\rho(X)=1$, and $K_X^2 =16$. 
\end{itemize}
\end{enumerate}
\end{cor}

\begin{proof}
The assertions (1) and (2) follow from Proposition \ref{p index 3}. 
The assertion (3) holds by  Proposition \ref{p index 2 pic2} and Theorem \ref{t index 2}. 
\end{proof}

\section{Boundedness over F-finite fields}\label{s bdd}

\subsection{Case $\rho(X)=2$}

\begin{lem}\label{l quad torsor}
Let $X$ be a regular del Pezzo surface. 
Assume that 
there exist an ample Cartier divisor $A$ and a non-trivial $\alpha_{-A}$-torsor $f: Y \to X$ 
such that $A \equiv -K_X$ and 
$Y$ is a normal quadric surface on $\P^3_{k_Y}$ for some field $k_Y$. 
Then the following hold: 
\begin{enumerate}
\item $p=2$. 
\item $h^1(X, \MO_X)=1$. 
\item $A \sim -K_X$. 
\end{enumerate}
\end{lem}

\begin{proof}

We have $K_Y \sim f^*(K_X -(p-1)A) \equiv pf^*K_X$. 
As $Y$ is a normal quadric surface, we get $p=2$. 
Thus (1) holds. 
We obtain 
\[
8\cdot [k_Y:k] = 
K_Y^2  =  (2f^*K_X)^2 = 2^2 \cdot \deg f \cdot K_X^2 = 8\cdot K_X^2, 
\]
which implies $[k_Y:k] = K_X^2$. 
Since $f: Y \to X$ is an $\alpha_{-A}$-torsor, we get 
\[
\chi(\MO_Y) -\chi(\MO_X)
= \sum_{m=1}^{2-1} \chi(X, mA) =  \chi(\MO_X) + \frac{1}{2}A \cdot (A -K_X)
=\chi(\MO_X) + K_X^2. 
\]
By $\chi(Y, \MO_Y) = h^0(Y, \MO_Y) = [k_Y:k] =K_X^2$, we get $\chi(X, \MO_X)=0$, i.e.,  $h^1(X, \MO_X)  = 1$. 
Thus (2) holds.

Let us show (3). 
Suppose that $A \equiv -K_X$ and $A \not\sim -K_X$. 
It is enough to derive a contradiction. 
We have that $h^0(X, -A)=0$ and $h^2(X, -A) = h^0(X, K_X +A)=0$. 
Hence 
\[
-h^1(X, -A) = \chi(X, -A) = \chi(X, \MO_X) + \frac{1}{2} (-A) \cdot (-A -K_X)
= 0 + 0 =0, 
\]
which contradicts $H^1(X, -A) \neq 0$ (Remark \ref{r alpha H^1 correspondence}).
Thus (3) holds. 
\end{proof}

\begin{prop}\label{p rho=2 h^1=1}
Let $X$ be a regular del Pezzo surface with $\rho(X) \geq 2$. 
Assume that there is an ample Cartier divisor $A$ such that $H^1(X, -A) \neq 0$. 
Then the following hold: 
\begin{enumerate}
\item $p=2$ and $\rho(X)=2$. 
\item $A \sim -K_X$. 
\item $h^1(X, \MO_X)=1$. 
\item $(X \times_k \ol{k})_{\red}^N \simeq \P^1 \times \P^1$. 
\item There exists a non-trivial $\alpha_{\omega_X}$-torsor $f: Y \to X$. 
\end{enumerate}
\end{prop}

\begin{proof}
By $H^1(X, -A) \neq 0$, $X$ is not geometrically normal 
(Proposition \ref{p dP basic2}). 
In particular, $p \in \{2, 3\}$ (Proposition \ref{p dP basic1}) and 
$\rho(X) =2$ \cite[Theorem 4.6]{Tan24}. 
Thus (1) holds. 

Take the largest integer $m>0$ such that $H^1(X, -mA) \neq 0$ and 
$H^1(X, -m'A)=0$ for every integer $m' > m$. 
In particular, $H^1(X, -mpA)=0$. 
Then we get an $\alpha_{-mA}$-torsor $f : Y \to X$ and hence 
\[
K_Y \sim f^*(K_X - m(p-1)A). 
\]
For the normalisation $\nu: Z \to Y$ and $g: Z \xrightarrow{\nu} Y \xrightarrow{f} X$, we get 
\[
K_Z + C \sim \nu^*K_Y \sim g^*(K_X - m(p-1)A). 
\]
Note that $Z$ is a projective normal $\Q$-factorial surface.

We may assume that $k$ is separably closed. 
By the classification of $(X \times_k \overline k)_{\red}^N$ \cite[Theorem 4.6]{Tan24}, 
we can find a contraction $\pi : Z \to B$ to a regular projective curve $B$. 
For the generic fibre $F$ of $\pi$, 
consider the restriction to $F$ of the following divisor: 
\[
-K_Z \sim   C  +g^*(-K_X + m(p-1)A). 
\]
Then $m=1$, $p=2$, and $C$ is $\pi$-vertical. 
Since $\rho(Z)=\rho(X)=2$,  $C$ is nef. 
As $g: Z \to X$ is an $\alpha_{M}$-torsor for some $M$, $Z$ is lci, and hence Gorenstein. 
Therefore, we have $-K_Z = A_Z + A'_Z$ for the ample Cartier divisors $A_Z := g^*(-K_X)$ and $A'_Z := -K_Z -A_Z \equiv C +g^*A$. 
Therefore, $Z \simeq  \P^1 \times \P^1$ (Proposition \ref{p index 2 pic2}). 
Then it holds that $-g^*K_X \sim \MO(1, 1) \sim g^*A$, $C=0$, and 
$(X \times_k \overline k)^N_{\red} \simeq \P^1 \times \P^1$. 
Then the assertion follows from Lemma \ref{l quad torsor}. 
\qedhere

\end{proof}

\subsection{Case $p=3$}

\begin{lem}\label{l Fujita pseff}
Let $Z$ be a projective normal $\Q$-factorial surface. 
Let $A$ be a nef and big Cartier divisor. 
Then $K_Z + 3A$ is pseudo-effective. 
\end{lem}

\begin{proof}
Taking the minimal resolution and the base change to the separable closure of $k$, we may assume that $Z$ is regular and $k$ is separably closed. 
By \cite[Theorem 1.1]{Tan18b}, we may assume that $k$ is algebraically closed. 
By taking the minimal resolution again, we may assume that $Z$ is smooth.

Suppose that $K_Z+3A$ is not pseudo-effective. 
We run a $(K_Z+3A)$-MMP: 
\[
Z =:Z_0 \to Z_1 \to \cdots \to Z_{\ell} =:W, 
\]
where $W$ denotes the end result. 
As $K_Z +3A$ is not pseudo-effective, 
we have a $(K_W+3A_W)$-negative Mori fibre space $\pi : W \to B$ for 
the push-forward $A_W$ of $A$. 
Note that the above MMP is a $K_Z$-MMP and $\pi : W \to B$ is a $K_W$-Mori fibre space. 
In particular,  each $Z_i$ is a smooth projective surface. 
However, we get a contradiction by taking an intersection with a line $\ell$ on $W=\P^2$ (resp.\ a fibre $\ell=\P^1$ of $\pi$) 
when $\dim B =0$ (resp.\ $\dim B=1$). 
\end{proof}

\begin{prop}\label{p easy vanish}
Let $X$ be a regular del Pezzo surface. 
Then $H^1(X, -A) = H^1(X, K_X+A)=0$ for 
every Cartier divisor $A$ satisfying $A \equiv -nK_X$ for some integer $n \geq 3$. 
\end{prop}

\begin{proof}
By Serre duality, it suffices to show $H^1(X, -A)=0$. 
Suppose $H^1(X, -A) \neq 0$. 
Replacing $A$ by $p^eA$ for some $e>0$ if necessary, 
we may assume that $H^1(X, -pA) =0$. 
Then we obtain an $\alpha_{-A}$-torsor $f: X' \to X$ satisfying 
$K_{X'} \sim f^*(K_X- (p-1)A) \equiv ( 1+ n(p-1)) f^*K_X$. 
For the normalisation $\nu :X'' \to X'$ and its conductor $C$, we obtain 
\[
K_{X''} +C \equiv \nu^*K_{X'} \equiv (1+n(p-1))\nu^*f^*K_X. 
\]
By $1+n(p-1) \geq 1+3 \cdot 1 =4$, this contradicts Lemma \ref{l Fujita pseff}. 
\end{proof}

\begin{prop}\label{p p=3 bdd}
Assume $p=3$. 
Let $X$ be a regular del Pezzo surface with $H^1(X, K_X) \neq 0$. 
Then 
there exists a non-trivial $\alpha_{K_X}$-torsor and 
the following hold for every non-trivial $\alpha_{K_X}$-torsor: 
\begin{enumerate}
\item $Y \simeq \P^2_{k_Y}$ for some field $k_Y$. 
\item $[k_Y:k] = 3K_X^2$. 
\item $K_X^2 = -3\chi(X, \MO_X)$ (equivalently, 
$h^1(X, \MO_X) = h^1(X, K_X) = 1+ \frac{K_X^2}{3})$. 
\end{enumerate}
\end{prop}


\begin{proof}
The existence of a non-trivial $\alpha_{K_X}$-torsor is ensured by $H^1(X, 3K_X)=0$. 
For the normalisation $\nu:Z \to Y$ of $Y$ and its conductor $C$, we have 
\[
K_Z + C = \nu^*K_Y \sim 3\nu^*f^*K_X. 
\]
By Lemma \ref{l Fujita pseff}, we get $C=0$, i.e., $Z = Y$. 
Then we obtain $K_Y \sim 3f^*K_X$, which implies $Y \simeq \P^2_{k_Y}$ (Corollary \ref{c r>1}). Thus (1) holds. 
Moreover, 
\[
9 [k_Y:k] = K_Y^2 = (3f^*K_X)^2 = 9 \cdot (\deg f) \cdot K_X^2 = 9 \cdot 3 \cdot K_X^2, 
\]
which implies $[k_Y:k]= 3K_X^2$. Thus (2) holds. 
On the other hand, we have 
\[
3K_X^2 = [k_Y : k] = 
\chi(Y, \MO_Y) = \sum_{i=0}^{p-1} \chi(X, -iK_X) = 3\chi(X, \MO_X) + 4K_X^2. 
\]
Thus (3) holds.  
\end{proof}

\subsection{Case $p=2$}

\begin{prop}\label{p p=2 2K-torsor}
Let $X$ be a regular del Pezzo surface. 
Assume that there exist a Cartier divisor $A$ 
satisfying 
$A \equiv -2K_X$ and $H^1(X, -A) \neq 0$. 
Then 
there exists a non-trivial $\alpha_{-A}$-torsor and 
the following hold for every non-trivial $\alpha_{-A}$-torsor $f: Y \to X$: 
\begin{enumerate}
\item $p=2$. 
\item $Y \simeq \P^2_{k_Y}$ for some field $k_Y$. 
\item $[k_Y:k] = 2K_X^2$. 
\item $K_X^2 = -2\chi(X, \MO_X)$. 
\item $h^1(X, -A) =1$. 
\end{enumerate}
\end{prop}


\begin{proof}
The existence of a non-trivial $\alpha_{-A}$-torsor is ensured by 
$H^1(X, -pA)=0$ (Remark \ref{r alpha H^1 correspondence}, Proposition \ref{p easy vanish}). 
By a similar argument as in Proposition \ref{p p=3 bdd}, 
we see that (1) and (2) hold. 
Moreover, $K_Y \sim f^*(K_X  - A) \equiv 3f^*K_X$. 
Hence 
\[
9 [k_Y:k] = K_Y^2 = (3f^*K_X)^2 = 9 \cdot (\deg f) \cdot K_X^2 = 9 \cdot 2 \cdot K_X^2, 
\]
which implies $[k_Y:k]= 2K_X^2$. Thus (3) holds. 
On the other hand, we have 
\[
\chi(Y, \MO_Y) - \chi(X, \MO_X) = \chi(X, A) = \chi(X, \MO_X) + 3K_X^2. 
\]
Hence 
\[
2K_X^2 -2\chi(X, \MO_X) =[k_Y:k] -2\chi(X, \MO_X)= \chi(Y, \MO_Y) -2\chi(X, \MO_X) =3K_X^2,
\]
which implies $-2 \chi(X, \MO_X) = K_X^2$. Thus (4) holds. 

Let us show (5). We have an exact sequence 
\[
0 \to \MO_X\to f_*\MO_Y \to \MO_X(A) \to 0. 
\]
Tensoring with $\MO_X(-A)$ and taking cohomologies, we obtain 
\[
0 = H^0(Y, \MO_Y(-f^*A)) \to H^0(X, \MO_X) 
\to H^1(X, -A) \to H^1(Y, \MO_Y(-f^*A)) =0, 
\]
which implies $h^1(X, -A) = h^0(X, \MO_X)=1$. 
Thus (5) holds. 
\end{proof}

\begin{prop}\label{p p=2 K-torsor}
Let $X$ be a regular del Pezzo surface with $\Pic X \equiv \Z K_X$. 
Assume that $p=2$, $H^1(X, K_X) \neq 0$, and $H^1(X, M)=0$ 
for every Cartier divisor $M$ satisfying $M \equiv 2K_X$. 
Then there exists a non-trivial $\alpha_{\omega_X}$-torsor 
and the following hold for every non-trivial $\alpha_{\omega_X}$-torsor $f: Y \to X$:   
\begin{enumerate}
\item $K_Y \sim 2f^*K_X$. 
\item $Y$ is a normal surface.  
\item[(3)] 
If $K_Y^2 = 8 [k_Y:k]$, then $K_X^2 = [k_Y:k]$ and $h^1(X, \MO_X) =1$. 
\item[(4)] 
If $K_Y^2 = 16 [k_Y:k]$, 
then $K_X^2 = 2[k_Y:k]$ and $h^1(X, \MO_X) =1 + \frac{1}{4}K_X^2$. 
\end{enumerate}
\end{prop}

\begin{proof}
By $H^1(X, K_X) \neq 0$ and $H^1(X, 2K_X)=0$, there exists a non-trivial $\alpha_{\omega_X}$-torsor. 
Take a non-trivial $\alpha_{\omega_X}$-torsor $f: Y \to X$. 
It holds that $K_Y \sim f^*(K_X + (p-1)K_X) = 2f^*K_X$. 
Thus (1) holds. 

Let us show (2). 
Suppose that $Y$ is not normal. 
For the normalisation $\nu : Z \to Y$, we have the induced morphisms: 
\[
g: Z \xrightarrow{\nu} Y \xrightarrow{f} X. 
\]
Since $Y$ is not normal, there exists a nonzero effective Weil divisor $C$ on $Z$ satisfying 
\[
K_Z + C \sim \nu^*K_Y \sim 2 \nu^*f^*K_X, 
\]
and hence 
\[
-K_Z \sim -2g^*K_X + C. 
\]
By Proposition \ref{p index 3}, we get $Z = \P^2_{k_Z}$ and $\MO_{\P^2}(1) \sim -g^*K_X$. 
Since $g: Z \to X$ is a double cover, 
$g$ is an $\alpha_{L}$-torsor for some Cartier divisor $L$ on $X$. 
By $K_Z \sim g^*(K_X+L)$ and $\MO_{\P^2}(1) \sim -g^*K_X$, 
we get $g^*L \equiv 2g^*K_X$, which implies $L \equiv 2K_X$. 
We then get $H^1(X, L) \neq 0$, which contradicts our assumption. 
Thus (2) holds. 

We have $K_Y^2 = 8\lambda [k_Y:k]$ for some $\lambda \in \{1, 2\}$ (Corollary \ref{c r>1}). 
Hence 
\[
8\lambda  [k_Y:k] = 
K_Y^2  \overset{{\rm (1)}}{=}  (2f^*K_X)^2 = 2^2 \cdot \deg f \cdot K_X^2 = 8\cdot K_X^2. 
\]
Thus $K_X^2 = \lambda [k_Y:k]$. 
Since $f: Y \to X$ is an $\alpha_{\omega_X}$-torsor, we get 
\[
\chi(Y, \MO_Y) -2\chi(X, \MO_X)
= \sum_{m=1}^{2-1} \frac{K_X^2}{2}m(m+1) = K_X^2. 
\]
By $\chi(Y, \MO_Y) = h^0(Y, \MO_Y) = [k_Y:k]$ (Theorem \ref{t index 2}), we get 
\[
[k_Y:k] -2\chi(X, \MO_X) 
= K_X^2  = \lambda    [k_Y:k].  
\]
Hence 
\[
-\chi(X, \MO_X) = \frac{\lambda-1}{2} [k_Y:k], 
\]
as required. 
\qedhere


\end{proof}

\subsection{Structure result and boundedness over F-finite fields}

We start by establishing the following structure result. 

\begin{thm}\label{t structure}
Let $X$ be a regular del Pezzo surface satisfying $H^1(X, \MO_X) \neq 0$. 
Then exactly one of the following four cases holds: 
\[
{\setlength{\arraycolsep}{3pt}
\begin{array}{c|c|c|c|c|c}
p
& \rho(X)
& \displaystyle-\frac{\chi(X,\MO_X)}{K_X^2}
& K_X^2
& (X \times_k \overline{k})_{\red}^N
& \begin{array}{c}\text{Torsor description}\end{array}\\
\hline
2 & 1, 2 & 0
& [k_Y:k]
& \begin{array}{c}
\P^2,\ \P^1\times\P^1,\\
\P(1,1,2)
\end{array}
& \begin{array}{c}
\exists\,\alpha_{\omega_X}\text{-torsor}\,\,Y\to X,\\
Y\text{:normal quadric}
\end{array}\\
\hline
2 & 1 & 1/4
& 2[k_Y:k]
& \P^2
& \begin{array}{c}
\exists\,\alpha_{\omega_X}\text{-torsor}\,\,Y\to X,\\
Y\text{:normal},\ K_Y^2=16[k_Y:k]
\end{array}\\
\hline
2 & 1 & 1/2
& \frac{1}{2}[k_Y:k] 
& \P^2
& \begin{array}{c}
\exists\,A\equiv-2K_X,\\
\exists\,\alpha_{-A}\text{-torsor}\,\,\P^2_{k_Y}\to X
\end{array}\\
\hline
3 & 1 & 1/3
& \frac{1}{3}[k_Y:k]
& \P^2
& \begin{array}{c}
\exists\,\alpha_{\omega_X}\text{-torsor}\,\,\P^2_{k_Y}\to X
\end{array}\\
\hline
\end{array}}
\]
Moreover, the following hold: 
\begin{enumerate}
\item $r_X=1$. In particular, if $\rho(X)=1$, then $\Pic X \equiv \Z K_X$. 
\item Every effective divisor on $X$ is nef. 
\end{enumerate}
\end{thm}

\begin{proof}
The assertion on the column \lq\lq $K_X^2$" (resp.\ \lq\lq $(X\times_k\overline{k})_{\red}^N$") 
follows from the column \lq\lq Torsor description" 
(resp.\ the other columns and \cite[Theorem 4.6]{Tan24}). 
Since $H^1(X, \MO_X)\neq 0$, 
Propositions \ref{p dP basic1} and
\ref{p dP basic2} imply $p\in\{2,3\}$.


The assertion (1) holds by Corollary \ref{c r>1}. 
If $\rho(X) \geq 2$, then the first case holds by Proposition \ref{p rho=2 h^1=1}. In particular, (2) holds. 
If $p=3$, then the assertion holds by 
Proposition \ref{p p=3 bdd}. 
In what follows, we assume $p=2$ and $\Pic X \equiv \Z K_X$.

If $H^1(X, -A) \neq 0$ for some Cartier $A$ with $A \equiv -2K_X$, 
then $K_X^2 = -2\chi(X, \MO_X)$ and the third case holds 
(Proposition \ref{p p=2 2K-torsor}). 
If $H^1(X, -A) = 0$ for every Cartier $A$ with $A \equiv -2K_X$, 
then the first or second case holds by Proposition \ref{p p=2 K-torsor}.
\qedhere


\end{proof}

\begin{thm}\label{t non-vanishing}
Let $X$ be a regular del Pezzo surface. 
Let $A$ be an ample Cartier divisor satisfying $H^1(X, -A) \neq 0$. 
Then exactly one of the following holds: 
\begin{enumerate}
\item $p=2$, $2A \sim -K_X$, $K_X^2 =16$, $H^1(X, \MO_X)=0$, and $h^1(X, -A)=1$. 
\item $p\in \{2, 3\}$, $A \sim -K_X$ and $h^1(X, -A) = h^1(X, \MO_X)>0$. 
\item $p\in \{2, 3\}$, $A \equiv -K_X$, $A \not\sim -K_X$, and $h^1(X, -A) = -\chi(X, \MO_X) >0$. 
\item $p=2$, $A \equiv -2K_X$, $-\chi(X, \MO_X)/K_X^2 = 1/2$, and $h^1(X, -A) =1$. 
\end{enumerate}
\end{thm}

\begin{proof}
Let $r_X$ be the index of $X$. 
If $r_X \geq 2$, then the assertion follows from Corollary \ref{c r>1}. 
In what follows, we assume $r_X=1$, and hence (1) no longer occurs. 
If $\rho(X) \geq 2$, then the assertion holds by Proposition \ref{p rho=2 h^1=1}. 

Let us treat the case when $\Pic X \equiv \Z K_X$. 
We then get $A \equiv - nK_X$ for some integer $n>0$. 
If $n \geq 3$, then $H^1(X, -A)=0$ by Proposition \ref{p easy vanish}. 
If $n=2$ (i.e., $A \equiv -2K_X$), then the assertion holds by 
Proposition \ref{p p=2 2K-torsor}. 
Hence we may assume $n=1$, i.e., $A \equiv -K_X$. 

If $A \sim -K_X$, then (2) holds. 
Assume  $A \not\sim -K_X$. 
We then get $H^0(X, -A)=0$ and $h^2(X, -A) = h^0(X, K_X+A) =0$, 
which implies 
\[
-h^1(X, -A) = \chi(X, -A) = \chi(X, \MO_X) + 
\frac{1}{2} (-A) \cdot (-A -K_X) = \chi(X, \MO_X).
\]
Therefore, $0< h^1(X, -A) = -\chi(X, \MO_X)$, and hence (3) holds.  
\end{proof}

\begin{thm}\label{t main}
Fix an integer  $v_0>0$. 
Then regular del Pezzo surfaces $X$ satisfying $K_X^2 \leq v_0$ are bounded. 

Specifically, there exists a flat projective morphism $\rho:V \to S$ 
of quasi-projective $\Z$-schemes 
which satisfies the following property: 
if $k$ is a field of characteristic $p>0$ and 
$X$ is a regular del Pezzo surface over $k$ such that $K_X^2 \leq v_0$, 
then there exists a cartesian diagram of schemes: 
\[
\begin{CD}
X @>>> V\\
@VV\alpha V @VV\rho V\\
\Spec\,k @>>> S, 
\end{CD}
\]
where $\alpha$ denotes the induced morphism. 

\end{thm}

\begin{proof}
Let $X$ be a regular del Pezzo surface. 
Recall that $-12K_X$ is very ample  \cite[Theorem 1.1]{Tan24}. 
We have 
\[
h^0(X, -12mK_X) = \chi(X, -12mK_X) 
= 1-h^1(X, \MO_X)+6m(12m+1) K_X^2. 
\]
In particular, $h^0(X, -12K_X) = 1-h^1(X, \MO_X) + 78 K_X^2 \leq 1 + 78 v_0$. 
Thus it is enough to find a uniform bound of $h^1(X, \MO_X)$.

We may assume $H^1(X, \MO_X) \neq 0$. 
Then the assertion follows from 
\[
h^1(X, \MO_X) -1 = -\chi(X, \MO_X) 
\leq \frac{1}{2} K_X^2 \overset{(\star)}{\leq} \frac{1}{2} v_0,
\]
where $(\star)$ is ensured by Theorem \ref{t structure}. 
\end{proof}

\begin{cor}\label{c F-fin bdd main}
Fix an integer $e_0>0$. 
Then regular del Pezzo surfaces $X$ over fields $k$ satisfying 
$[k:k^p] \leq p^{e_0}$  are bounded. 
\end{cor}

\begin{proof}
The assertion follows from Theorem \ref{t main} and \cite[Theorem 1.2(1)]{Tan24}.
\end{proof}

\section{Examples with small volumes}\label{s examples}

\subsection{Construction}

\begin{prop}\label{p common construction}
Let $k$ be a field of characteristic $p>0$. 
Take a purely inseparable field extension
$k
\subset
k':=k(\alpha_1,\ldots,\alpha_m)$ 
such that $[k':k]=p^m$ 
and $\alpha_i^p=:a_i\in k$ for every 
$1 \leq i\leq m$. 
Let $V$ be a regular projective variety over $k$ such that 
$H^0(V, \MO_V)=k$ and $Y:=V_{k'}$ is regular. 
Fix a Cartier divisor $L_V$ on $V$ and linearly independent elements
$q_1,\ldots,q_m
\in
H^0(V,\mathcal O_V(L_V))$ 
having no common zero.
For each $i$, let
\[
\partial_i: \MO_Y =\MO_V \otimes_k k' \to \MO_V \otimes_k k'=\MO_Y
\]
be the $\MO_V$-linear derivation satisfying $\partial_i(\alpha_j)=\delta_{ij}$. 
For the pullback $L$ of $L_V$ on $Y$, set 
\[
D:=
\sum_{i=1}^m q_i\partial_i:
\mathcal O_Y
\to
\mathcal O_Y(L), 
\]
which is an $\MO_V$-linear derivation. 
Let $f : Y \to Y^D =:X$ be the induced morphism. 
Then the following hold:



\begin{enumerate}
\item
We have finite surjective morphisms 
\[
Y = V_{k'} \xrightarrow{f}X\xrightarrow{g}V
\]
of regular projective varieties such that 
the composition $V_{k'} \to V$ is the base-change morphism induced by $k'/k$ and 
$f$ is a non-trivial $\alpha_{-M}$-torsor for $M:=g^*L_V$. 
\item $H^0(X,\mathcal O_X)=k$. 
\item $\omega_Y
\simeq
f^*
\left(
\omega_X\otimes\mathcal O_X(-(p-1)M)
\right).$ 
\item
There exists a filtration $\mathcal O_X=:\mathcal E_0
\subset
\mathcal E_1
\subset
\cdots
\subset
\mathcal E_{p-1}:=f_*\mathcal O_Y
$ such that
$\mathcal E_j/\mathcal E_{j-1}
\simeq
\mathcal O_X(jM)$. 

\item
$p^m\chi(V,\mathcal O_V)
=
\sum_{j=0}^{p-1}
\chi(X,\mathcal O_X(jM))$. 
\end{enumerate}
\end{prop}

\begin{proof}
Let us show (1). 
As $D$ is $\MO_V$-linear, we get $D(\MO_V)=0$, 
which induces the required factorisation $Y = V_{k'} \xrightarrow{f} X \xrightarrow{g} V$. 
We have $L = f^*g^*L_V = f^*M$ and $D^p=0$. 
Since $q_1, ..., q_m$ have no common zero, 
$D$ is nowhere vanishing. 
Hence $f$ is an $\alpha_{-M}$-torsor by Proposition \ref{p derivation to torsor}. 
In particular, $f$ is faithfully flat, and hence $X$ is a regular projective variety. 
Thus (1) holds.




Let us show (2). Recall that $H^0(Y,\mathcal O_Y)=k'$. 
Take $\lambda\in H^0(X, \MO_X)$. 
Then $\lambda \in k'$ and $D(\lambda)=0$. 
Hence 
\[
0 = D(\lambda) = \sum_{i=1}^m \partial_i(\lambda)q_i \qquad \text{in}\qquad 
H^0(Y, L) = H^0(V, L_V) \otimes_k k'.
\]
Note that $\partial_i(\lambda) \in k'$ for every $i$. 
Since $q_1,\ldots,q_m$ are linearly independent over $k'$, we obtain
$\partial_i(\lambda)=0$ for every $i$. 
Since $\alpha_1,\ldots,\alpha_m$ form a $p$-basis of
$k'/k$, this implies $\lambda\in k$. Thus (2) holds.

By using (1), the assertion (3) and (4) follow from 
Proposition \ref{p torsor} and \cite[the proof of Proposition 1.7]{Eke88}, respectively. 
Then (4) and $\chi(Y, \MO_Y) =[k':k]\chi(V, \MO_V) =p^m\chi(V, \MO_V)$ imply (5). 
\qedhere


\end{proof}

\begin{prop}\label{p H1 dimension common construction}
We use the same notation and assumptions of Proposition
\ref{p common construction}. Assume that
$H^0(V, \MO_V(-L_V))=H^1(V,\MO_V(-L_V))
=
0.$ 
Then
\[
h^1(X,\MO_X(-M))
\geq
\binom{m+p-2}{p-2}.
\]
If the natural map $\Sym^{p-1}\left(
H^0(V,\MO_V(L_V))\right)
\to
H^0(V,\MO_V((p-1)L_V))$ is injective, then
\[
h^1(X,\MO_X(-M))
=
\binom{m+p-2}{p-2}.
\]
\end{prop}

\begin{proof}
Recall that
$D=\sum_{i=1}^m q_i\partial_i:
\MO_Y\to\MO_Y(L)$. 

\setcounter{step}{0}

\begin{step}\label{s1 H1 dimension common construction}
We have the following exact sequence of $\MO_X$-modules:
\[
0
\to
\MO_X(-M)
\to
f_*\MO_Y(-L)
\xrightarrow{D}
f_*\MO_Y
\xrightarrow{D^{p-1}}
f_*\MO_Y((p-1)L).
\]
\end{step}

\begin{proof}[Proof of Step \ref{s1 H1 dimension common construction}]
Note that the pushforward 
$D:
f_*\MO_Y
\to
f_*\MO_Y(L)$
is an $\MO_X$-linear homomorphism, 
which induces another $\MO_X$-module homomorphism 
$D: f_*\MO_Y(f^*N) \to f_*\MO_Y(L+f^*N)$ for every Cartier divisor $N$ on $X$. 
Then we obtain $\Ker(D: f_*\MO_Y(-L)
\to 
f_*\MO_Y) = \MO_X(-M)$. 

It suffices to verify exactness at $f_*\MO_Y$. 
The problem is local on $X$. 
By Lemma \ref{l additive slice} and Lemma \ref{l unit slice}, 
we may assume that 
\[
\MO_Y
=
\bigoplus_{j=0}^{p-1}\MO_Xt^j,
\qquad
D(t)=1
\]
for some $t \in \MO_Y$. We obtain 
$\Ker D^{p-1} 
=
\bigoplus_{j=0}^{p-2}\MO_Xt^j
=
\Image D.$ 
This completes the proof of Step \ref{s1 H1 dimension common construction}. 
\end{proof}

\begin{step}\label{s2 H1 dimension common construction}
$H^1(X, \MO_X(-M)) \simeq \Ker( D^{p-1} : H^0(Y, \MO_Y) \to H^0(Y, \MO_Y((p-1)L)))$. 
\end{step}

\begin{proof}[Proof of Step \ref{s2 H1 dimension common construction}]
For the coherent $\MO_X$-module $\mathcal C
:=
\Image\left(
D:f_*\MO_Y(-L)\to f_*\MO_Y
\right)$, we have the following exact sequences: 
\[
0
\to
\MO_X(-M)
\to
f_*\MO_Y(-L)
\to
\mathcal C
\to
0
\]
\[
0 \to \cC \to f_*\MO_Y
\xrightarrow{D^{p-1}}
f_*\MO_Y((p-1)L). 
\]
These exact sequences, together with $H^0(Y,\MO_Y(-L))
=
H^1(Y,\MO_Y(-L))
=
0$, imply 
\[
H^1(X,\MO_X(-M))
\simeq
H^0(X,\mathcal C) \simeq \Ker\left(
D^{p-1}:H^0(Y, \MO_Y) 
\to
H^0(Y, \MO_Y((p-1)L))
\right). 
\]
This completes the proof of Step \ref{s2 H1 dimension common construction}. 
\end{proof}

\begin{step}\label{s3 H1 dimension common construction}
$h^1(X,\MO_X(-M))
\geq
\binom{m+p-2}{p-2}$.
\end{step}

\begin{proof}[Proof of Step \ref{s3 H1 dimension common construction}]
For a multi-index $I=(i_1,\ldots,i_m)\in\mathbb Z_{\geq0}^m$, 
set
$|I|:=i_1+\cdots+i_m$ and  $\alpha^I:=\alpha_1^{i_1}\cdots\alpha_m^{i_m}.$ 
Every element of $\bigoplus_{|I|\leq p-2}k\alpha^I$ 
is annihilated by $D^{p-1} : k' = H^0(Y, \MO_Y) \to H^0(Y, \MO_Y((p-1)L))$. 
Consequently,
\[
h^1(X,\MO_X(-M)) = \dim_k \Ker(D^{p-1}: k' \to H^0(Y, \MO_Y((p-1)L))) 
\]
\[
\geq 
\dim_k
\left(
\bigoplus_{|I|\leq p-2}k\alpha^I
\right)
=
\binom{m+p-2}{p-2}.
\]
This completes the proof of Step \ref{s3 H1 dimension common construction}. 
\end{proof}

\begin{step}\label{s4 H1 dimension common construction}
If $\Sym^{p-1}
(H^0(V,\MO_V(L_V)))
\to
H^0(V,\MO_V((p-1)L_V))$ is injective, then 
\[
h^1(X,\MO_X(-M))
=
\binom{m+p-2}{p-2}.
\]
\end{step}

\begin{proof}[Proof of Step \ref{s4 H1 dimension common construction}]
Assume  that
$\Sym^{p-1}(
H^0(V,\MO_V(L_V)))
\to
H^0(V,\MO_V((p-1)L_V))$ is injective. 
Then the composition 
\[
\Sym^{p-1}
(kq_1 \oplus \cdots \oplus k q_m) \hookrightarrow \Sym^{p-1}
(H^0(V,\MO_V(L_V)))
\to H^0(V,\MO_V((p-1)L_V))
\]
is also injective. Hence the elements
$q^I:=q_1^{i_1}\cdots q_m^{i_m} \in H^0(V,\MO_V((p-1)L_V))$ 
with $|I|=p-1$ are linearly independent over $k$.
We have
\[
D^{p-1}
=
\sum_{|I|=p-1}
\binom{p-1}{I}q^I\partial^I, \qquad \text{where}\qquad 
\partial^I
:=
\partial_1^{i_1}\cdots\partial_m^{i_m}.
\]
Take $\lambda \in \Ker D^{p-1}$. 
As $\binom{p-1}{I}\neq0$ in $k$ for every $I$ satisfying $|I|=p-1$, the equality 
$D^{p-1}(\lambda)=0$ implies
\[
\partial^I(\lambda)=0
\qquad
\text{for every }I\text{ satisfying }|I|=p-1.
\]
We can write 
\[
\lambda  = \sum_{0 \leq i_j \leq p-1} c_{i_1, ..., i_m} \alpha_1^{i_1}\cdots\alpha_m^{i_m} \in \bigoplus_{0 \leq i_j \leq p-1} k\alpha_1^{i_1}\cdots\alpha_m^{i_m}  =k', 
\]
and hence we obtain $c_{i_1, ..., i_m}=0$ when 
$|I| = i_1 + \cdots + i_m \geq p-1$. 
Therefore, 
$\Ker D^{p-1}
=
\bigoplus_{|I|\leq p-2}k\alpha^I$. 
By the same computation as in the previous step, we obtain 
$h^1(X,\MO_X(-M))
=
\binom{m+p-2}{p-2}$. 
This completes the proof of Step \ref{s4 H1 dimension common construction}. 
\end{proof}
Step \ref{s3 H1 dimension common construction} and Step \ref{s4 H1 dimension common construction} complete the proof of 
Proposition \ref{p H1 dimension common construction}. 
\end{proof}

\begin{prop}\label{p common construction2}
Let $k$ be a field of characteristic $p>0$. 
Take a purely inseparable field extension
\[
k
\subset
k':=k(\alpha_1,\ldots,\alpha_m)
\]
such that $[k':k]=p^m$ and $\alpha_i^p=:a_i\in k$ for every 
$1 \leq i\leq m$. 
Let $V$ be an $n$-dimensional regular projective variety over $k$ such that 
$H^0(V, \MO_V)=k$ and $Y :=V_{k'}$ is regular.  
Fix a Cartier divisor $H$ on $V$ and integers $d>0$ and $r>0$. 
Assume that  
\begin{enumerate}
\renewcommand{\labelenumi}{(\alph{enumi})}
\item $-K_V \sim r H$, 
\item $|dH|$ is base point free, and 
\item $n < m \leq h^0(V, dH)$. 
\end{enumerate}
Then there exist finite surjective morphisms 
\[
Y =V_{k'} \xrightarrow{f} X \xrightarrow{g} V
\]
of $n$-dimensional regular projective varieties  
such that the composition $Y=V_{k'} \to V$ is the base change and the following properties hold for $A := g^*H$: 
\begin{enumerate}
\item $f$ is a non-trivial $\alpha_{-dA}$-torsor. 
Moreover, $\deg f = p$ and $\deg g = p^{m-1}$. 
In particular, $A^n = p^{m-1}H^n$. 
\item $H^0(X,\mathcal O_X)=k$. 
\item
There exists a filtration $\mathcal O_X=:\mathcal E_0
\subset
\mathcal E_1
\subset
\cdots
\subset
\mathcal E_{p-1}:=f_*\mathcal O_Y
$ such that
$\mathcal E_j/\mathcal E_{j-1}
\simeq
\mathcal O_X(jdA)$. 
\item $K_Y
\sim 
f^*
\left(
K_X -(p-1)dA
\right)$ and  $-K_X \equiv (r-(p-1)d)A$.

\item 
For $s := r-(p-1)d$, it holds that 
\[
p^m\chi(V,\mathcal O_V)
=
\sum_{j=0}^{p-1}
\chi(X,\mathcal O_X(jdA)).
\]
If $n=2$, then $p^m\chi(V, \MO_V) =   p\chi(X, \MO_X) + \frac{A^2}{2} \sum_{j=1}^{p-1} jd(jd+s)$. 
\item 
$H^1(X, \MO_X(-dA)) \neq 0$. 
If $H^1(V,\mathcal O_V(-dH))=0,$ 
then 
$h^1(X,\mathcal O_X(-dA))\geq 
\binom{m+p-2}{p-2}$.  
If $p=2$ and $H^1(V,\mathcal O_V(-dH))=0,$ 
then 
$h^1(X,\mathcal O_X(-dA))=1$. 
\end{enumerate}
\end{prop}

\begin{proof}
The assertions (1)-(3) follow from  Proposition \ref{p common construction} by setting $L_V := dH$ and $M := g^*L_V = dg^*H=dA$. 
Concerning (4), 
the linear equivalence 
$K_Y
\sim 
f^*
\left(
K_X -(p-1)dA
\right)$ follows from Proposition \ref{p common construction}, 
and it induces $-K_X \equiv (r-(p-1)d)A$ by (a). 
The assertion (5) holds by Proposition \ref{p common construction} and the Riemann-Roch theorem.

Let us show (6). 
We have $H^1(X, \MO_X(-dA)) \neq 0$ by (1) and 
Remark \ref{r alpha H^1 correspondence}. 
The second assertion follows from 
Proposition \ref{p H1 dimension common construction}.

Assume $p=2$ and $H^1(V,\mathcal O_V(-dH))=0$. 
By (3), we have an exact sequence
\[
0
\to
\mathcal O_X(-dA)
\to
f_*\mathcal O_Y(-dH)
\to
\mathcal O_X
\to
0.
\]
By $H^0(V, \MO_V(-dH)) = H^1(V,\mathcal O_V(-dH))=0$, 
we get an 
isomorphism
$H^0(X,\mathcal O_X)
\xrightarrow{\simeq}
H^1(X,\mathcal O_X(-dA))$. 
Thus (6) holds.
\qedhere

\end{proof}

\subsection{Examples}

\begin{ex}[index 2, $K_X^2=16$]\label{e index2 K^2=16}
Applying Proposition \ref{p common construction2} for 
\[
(V,H,p, r, d,m)
=
(\mathbb P^2_k,\mathcal O_{\mathbb P^2_k}(1),2,3, 1,3),
\]
we obtain a regular del Pezzo surface $X$ over 
a field $k$ of characteristic two such that $H^0(X, \MO_X)=k$, 
$-K_X \sim 2A$, $h^1(X, \MO_X(-A))=1$, and $K_X^2=16$ (note that 
$-K_X \equiv 2A$ implies $-K_X \sim 2A$ by 
Proposition \ref{p dP basic2}(2) and 
Theorem \ref{t index 2}(3)). 
Moreover, $|A|$ is base point free and $h^0(X, A) = 7$, because 
Theorem \ref{t index 2} implies 
\[
h^0(X, A) = \chi(X, A) = \chi(X, \MO_X) + \frac{1}{2}A \cdot (A -K_X) = 1 + \frac{1}{2} \cdot 3 \cdot 4 = 7. 
\]
\end{ex}

\begin{ex}[$\rho=2$, $H^1(\MO_X) \neq 0$]\label{e rho=2 irregular}

Applying Proposition \ref{p common construction2} for $m \in \{3, 4\}$ and 
\[
(V,H,p,r, d)
=
(\mathbb P^1_k\times_k\mathbb P^1_k,
\mathcal O_{\mathbb P^1_k\times_k\mathbb P^1_k}(1,1),2,2, 1),
\]
we obtain a regular del Pezzo surface $X$ over 
a field $k$ of characteristic two  such that $H^0(X, \MO_X)=k$, 
$\rho(X)=2, K_X^2=2^m,$ and $H^1(X, \MO_X) \overset{(\star)}{\neq} 0$. 
\end{ex}


\begin{proof}[Proof of $(\star)$]
By Proposition \ref{p common construction2}, we have $-K_X \equiv A$ and 
$H^1(X, \MO_X(-A)) \neq 0$. 
Suppose that $H^1(X, \MO_X)=0$. 
Then $-K_X \equiv A$ implies $-K_X \sim A$ 
(Proposition \ref{p dP basic2}(2)), 
which leads to the following contradiction: 
\[
0 \neq h^1(X, \MO_X(-A)) = h^1(X, \MO_X(K_X)) = h^1(X, \MO_X) =0. 
\]
Therefore, $H^1(X, \MO_X)\neq 0$. 
\end{proof}

\begin{prop}[$p=2$, pathological $\alpha_{\omega_X}$-torsor]\label{p common pathological}
Fix an integer $m$ satisfying $3 \leq m \leq 7$. 
Then there exists a regular del Pezzo surface $X$ over a field $k$ of characteristic two such that 
\begin{enumerate}
\item $H^0(X, \MO_X)=k$, 
\item $\Pic X \equiv \Z K_X$, 
\item $K_X^2 =2^{m+1}$, 
\item $h^1(X, \MO_X) = h^1(X, K_X) =2^{m-1}+1$, and 
\item $H^1(X, -M)=0$ for every Cartier divisor $M$ satisfying $M \equiv -2K_X$. 
\end{enumerate}
\end{prop}

\noindent 
Let $f : Y \to X$ be a non-trivial $\alpha_{\omega_X}$-torsor. 
Then Proposition \ref{p p=2 K-torsor} implies that 
$Y$ is a projective normal Gorenstein surface such that 
$K_Y^2 \geq 16$ and $K_Y \sim 2f^*K_X$ (hence of index $2$).

\begin{proof}
By Example \ref{e index2 K^2=16}, 
there exist a regular del Pezzo surface $V$ over a field $k$ 
of characteristic two and 
an ample Cartier divisor $H$ on $V$ such that $H^0(V, \MO_V)=k$, 
$-K_V \sim 2H$, $h^1(V, \MO_V(-H))=1$, $H^2 =4$, $h^0(V, H)=7$, and $|H|$ is base point free. 
Applying Proposition \ref{p common construction2} for 
$(V,H,p,r,d,m)$ by setting $r:=2$ and $d :=1$, 
we obtain a regular del Pezzo surface $X$ and an ample Cartier divisor $A$ on $X$ such that $-K_X \equiv A$, 
$K_X^2 = A^2 = p^{m-1}H^2 = 2^{m+1}$, and 
\[
2^m = p^m\chi(V, \MO_V) = p\chi(X, \MO_X) + \frac{A^2}{2} \cdot 1 \cdot 2 = 2(1-h^1(X, \MO_X)) + 2^{m+1}. 
\]
Hence $h^1(X, \MO_X)=2^{m-1}+1$. 
As $\rho(X) \leq \rho(V_{k'})=1$, we get $\rho(X)=1$. 
This, together with $H^1(X, \MO_X) \neq 0$, implies $\Pic X \equiv \Z K_X$ (Corollary \ref{c r>1}).

Suppose that $H^1(X, -M) \neq 0$ for some Cartier divisor $M$ satisfying $M \equiv -2K_X$. 
Then $K_X^2 =-2\chi(X, \MO_X)$ (Proposition \ref{p p=2 2K-torsor}), 
which contradicts $K_X^2 \overset{{\rm (3)}}{=} 2^{m+1} = 4 \cdot 2^{m-1}
 \overset{{\rm (4)}}{=}
-4\chi(X, \MO_X)$. 
\end{proof}

\begin{ex}[$p=2$, $H^1(X, 2K_X) \neq 0$]\label{e H1(2KX)}
Applying Proposition \ref{p common construction2} for $3 \leq m \leq 6$ and 
\[
(V,H,p,r, d)
=
(\mathbb P^2_k,\mathcal O_{\mathbb P^2_k}(1),2,3,2),
\]
we obtain a regular del Pezzo surface $X$ and an ample Cartier divisor $A$ such that $H^0(X, \MO_X)=k$, 
$\rho(X)=1, -K_X \equiv A, K_X^2=2^{m-1},$ 
\begin{enumerate}
\item[($\star$)] 
$h^1(X, 2K_X) =1$, and $h^1(X, \MO_X) = 2^{m-2}+1$. 
\end{enumerate}
In particular, by setting $m=6$, we obtain $K_X^2 = 32$ and $h^1(X, \MO_X) = 17$. 
\end{ex}

\begin{proof}[Proof of $(\star)$]
We have the factorisation $F: X \to Y =\P^2_{k'} \xrightarrow{f} X$ of the absolute Frobenius morphism of $X$, and hence 
the numerical equivalence $-K_X \equiv A$ implies the linear equivalence 
$-2K_X = -F^*K_X \sim F^*A =2A$. Therefore, we obtain 
$h^1(X, \MO_X(2K_X)) = h^1(X, \MO_X(-2A))=1$, where the latter equality is ensured by Proposition \ref{p common construction2}. 

Again by Proposition \ref{p common construction2}, 
we obtain 
\[
2^m=2^m\chi(V, \MO_V) = 2\chi(X, \MO_X) + \frac{A^2}{2} \cdot 2\cdot (2+1) = 2(1-h^1(X, \MO_X)) + 3 \cdot 2^{m-1}, 
\]
which implies $h^1(X, \MO_X) = 1 + 2^{m-2}$.     
\end{proof}

\begin{ex}[$p=3$, $h^1(X, \MO_X) =4$]\label{e p=3 irregular}
Applying Proposition \ref{p common construction2} for 
\[
(V,H,p,r, d,m)
=
(\mathbb P^2_k,\mathcal O_{\mathbb P^2_k}(1),3,3,1,3),
\]
we obtain a regular del Pezzo surface $X$ over a field $k$ of characteristic three such that 
\begin{enumerate}
\item[($\star$)] $K_X^2=9$ and $h^1(X,\mathcal O_X)=4.$ 
\end{enumerate}
\end{ex}

\begin{proof}[Proof of $(\star)$]
By Proposition \ref{p common construction2}, we have $-K_X \equiv A$ and $A^2 = p^2 H^2 = 9$ for some Cartier divisor $A$ on $X$. 
In particular, $K_X^2 = A^2 = 9$. 
Again by Proposition \ref{p common construction2}, we obtain  
\[
27=p^m \chi(V, \MO_V) = 3\chi(X, \MO_X) + \frac{A^2}{2}\sum_{j=1}^2 j(j+1) = 3(1-h^1(X, \MO_X)) + 36. 
\]
Hence $h^1(X, \MO_X)=4$. 
\end{proof}

The following example is a regular K3-like surface violating Kodaira vanishing. 
Note that Kodaira vanishing holds for geometrically integral K3-like surfaces \cite[Theorem 3.4]{FanoI}, which plays an essential role in classification of Fano threefolds in positive characteristic.

\begin{prop}[K3-like surfaces violating Kodaira vanishing]
\label{p K3 KV fail}
There exist a field $k$ of characteristic two, 
a regular projective surface $X$ over $k$, and an ample Cartier divisor $A$ such that 
$H^0(X,\MO_X)=k, K_X\sim0, H^1(X,\MO_X)=0, A^2=8,$ and
\[
h^1(X,\MO_X(A))
=
h^1(X,\MO_X(-A))
=
1.
\]
\end{prop}

\begin{proof}
Fix an algebraically closed field $\bF$ of characteristic two and 
a smooth del Pezzo surface $V_0$ over $\bF$ with $K_{V_0}^2 = 2$. 
For the purely transcendental extension $k:=\bF(a_1,a_2,a_3)$ of $\bF$, set 
$V:=V_0\times_{\bF}k$ and $H:=-K_V$. 
Applying Proposition \ref{p common construction2} for
$(p,r, d,m)=(2,1,1,3)$, 
we obtain a regular projective surface $X$ and an ample Cartier divisor $A$ satisfying the required properties except for 
\begin{equation}\label{e1 K3 KV fail}
K_X \sim 0 \qquad \text{and}\qquad H^1(X, \MO_X)=0. 
\end{equation}
It suffices to show (\ref{e1 K3 KV fail}). 

Proposition \ref{p common construction2} implies $K_X \equiv 0$ and 
\[
8
=
2^3\chi(V,\MO_V)
=
2\chi(X,\MO_X)
+
\frac{A^2}{2} = 2\chi(X,\MO_X)
+4. 
\]
Hence we get $\chi(X, \MO_X)=2$, and hence $h^0(X, K_X)=h^2(X, \MO_X)>0$. 
This, together with $K_X \equiv 0$, implies $K_X \sim 0$ and 
$h^2(X, \MO_X)=1$. 
We then obtain 
\[
2 = \chi(X, \MO_X) = h^0(X, \MO_X) -h^1(X, \MO_X) + h^2(X, \MO_X)=2-h^1(X, \MO_X). 
\]
Therefore, $H^1(X, \MO_X)=0$. 
\end{proof}

\section{Unbounded examples}\label{s unbounded}

\subsection{Case $p=3$}

\begin{thm}\label{t unbounded p=3}
Take an integer $n\geq 1$ and a field $k$ of characteristic three satisfying 
$[k:k^3] \geq 3^{n+3}$. 
Then there exists a regular del Pezzo surface
$X$ over $k$ such that $K_{X}^2=3^{n+2}$ 
and 
$h^1(X,\MO_{X})
=
h^1(X,K_{X})
=
1+3^{n+1}$.
\end{thm}

\begin{proof}
By $[k^{1/3}:k] =[k:k^3]\geq 3^{n+3}$, 
we can find an intermediate field 
\[
k
\subset
k' := k
\left(
\alpha_0,\alpha_1,\alpha_2,
\beta_1,\ldots,\beta_n
\right)
 \subset k^{1/3}
\]
such that $[k':k]=3^{n+3}$. 
Set $a_i := \alpha_i^3 \in k$ and $b_r :=\beta_r^3 \in k$
for every $0\leq i\leq 2$ and $1\leq r\leq n$.
Set
$Y:=\P^2_{k'}=\Proj k'[x_0,x_1,x_2]$.
For $0\leq r\leq n$, set
$k_r
:=
k
\left(
\alpha_0,\alpha_1,\alpha_2,
\beta_1,\ldots,\beta_r
\right)$.
For example,
\[
k_0=k(\alpha_0,\alpha_1,\alpha_2)\qquad \text{and}\qquad 
k_n=k
\left(
\alpha_0,\alpha_1,\alpha_2,
\beta_1,\ldots,\beta_n
\right)=k'.
\]
\setcounter{step}{0}

\begin{step}\label{s1 unbounded p=3}
There exists a $k$-derivation
\[
D:k'[x_0,x_1,x_2]
\to
k'[x_0,x_1,x_2]
\]
satisfying the following conditions. 
\begin{enumerate}
\item $D(x_0)=D(x_1)=D(x_2)=0$.
\item $D(\alpha_i)=x_i$ for every $0\leq i\leq 2$.
\item For every $1\leq r\leq n$, there exists $F_r\in k_{r-1}$ such that
$D(\beta_r)=F_rx_1$.
\item $F_1=\alpha_0$, i.e., $D(\beta_1)=\alpha_0x_1$.
\item If $1\leq r<n$, then
$F_{r+1}
=
-\beta_r^2F_r+F_r^3\alpha_1^2$. 
\end{enumerate}
\end{step}

\begin{proof}[Proof of Step \ref{s1 unbounded p=3}]
By (1), it is enough to find a $k$-derivation
$D:k'\to k'[x_0,x_1,x_2]$
satisfying (2)--(5). 
By $k_0=k(\alpha_0,\alpha_1,\alpha_2)$, 
there exists a unique $k$-derivation
$D_0:k_0\to
k_0x_0\oplus k_0x_1\oplus k_0x_2$
satisfying (2). There is a unique $k$-derivation 
$D_1:k_1\to
k_1x_0\oplus k_1x_1\oplus k_1x_2$
such that $D_1|_{k_0}=D_0$ and
$D_1(\beta_1)=\alpha_0x_1$.

Fix $1 \leq r < n$ and assume that 
a $k$-linear derivation $D_r: k_r \to k_rx_0\oplus k_rx_1\oplus k_rx_2$ has been defined. Since
$F_{r+1}
=
-\beta_r^2F_r+F_r^3\alpha_1^2
\in k_r$,
there exists a unique $k$-derivation 
$D_{r+1}:k_{r+1}\to
k_{r+1}x_0\oplus k_{r+1}x_1\oplus k_{r+1}x_2$
such that $D_{r+1}|_{k_r}=D_r$ and
$D_{r+1}(\beta_{r+1})=F_{r+1}x_1$.
We then obtain the  $k$-derivation  $D_n: k_n \to k_nx_0\oplus k_nx_1\oplus k_nx_2$, which is nothing but the required $k$-derivation by $k' = k_n$.
This completes the proof of Step \ref{s1 unbounded p=3}. 
\end{proof}

\begin{step}\label{s2 unbounded p=3}
It holds that $D^3=0$.
\end{step}

\begin{proof}[Proof of Step \ref{s2 unbounded p=3}]
We have $D^2(\alpha_i) \overset{{\rm (2)}}{=}D(x_i)\overset{{\rm (1)}}{=}0$
for every $0\leq i\leq 2$. 
As $k'=k_n$ and 
$D^3 :k'[x_0, x_1, x_2] \to k'[x_0, x_1, x_2]$ is a $k$-derivation, 
it suffices to show 
$D^3(k_n)=0$. 
We prove $D^3(k_r)=0$ by induction on $r$. 
The base case $r=1$ of this induction follows from 
$k_1 = k(\alpha_0, \alpha_1, \alpha_2, \beta_1)$ and 
$D^3(\beta_1) \overset{{\rm (4)}}{=} D^2(\alpha_0x_1) 
\overset{{\rm (1)}}{=}x_1 D^2(\alpha_0) = 0$. 


Fix $1\leq r<n$ and assume that $D^3(k_r)=0$. 
By $k_{r+1} = k_r(\beta_{r+1})$, it is enough to show 
$D^3(\beta_{r+1})=0$. We have 
\[
\begin{aligned}
D^3(\beta_{r+1})
&\overset{{\rm (3)}}{=}
D^2(F_{r+1}x_1) \\
&\overset{{\rm (5)}}{=}D^2(x_1(-\beta_r^2F_r+F_r^3\alpha_1^2))\\
&\overset{{\rm (3)}}{=}D^2(-\beta_r^2D(\beta_r)+F_r^3\alpha_1^2x_1)\\
&= -D^2(\beta_r^2 D(\beta_r)) + D^2(F_r^3x_1 \alpha_1^2). 
\end{aligned}
\]
The problem is reduced to showing 
\[
{\rm (i)}\,\,D^2(F_r^3x_1 \alpha_1^2)=2F_r^3x_1^3 \qquad \text{and}\qquad 
{\rm (ii)}\,\,D^2(\beta_r^2 D(\beta_r))=2F_r^3x_1^3.
\]
By $D(F_r^3)=D(x_1)=0$, we have 
\[
D^2(F_r^3x_1 \alpha_1^2) = F_r^3x_1 D^2(\alpha_1^2)=
F_r^3x_1 ( 2\alpha_1 D^2(\alpha_1) + 2D(\alpha_1)^2)
\overset{{\rm (2)}}{=}2F_r^3x_1^3.
\]
Thus (i) holds. 
Substituting $f=\beta_r$ and $g=D(\beta_r)$ 
into 
\[
D^2(f^2g)
=2D(f)^2g+2fD^2(f)g
+4fD(f)D(g)+f^2D^2(g),
\]
we obtain  
\[
D^2(\beta_r^2 D(\beta_r)) = 
2D(\beta_r)^2D(\beta_r) + 2\beta_rD^2(\beta_r) D(\beta_r)
+4\beta_r D(\beta_r) D^2(\beta_r) + \beta_r^2 D^3(\beta_r)
\]
\[
=2D(\beta_r)^3\overset{{\rm (3)}}{=}2 (F_rx_1)^3 = 2F_r^3x_1^3. 
\]
Thus (ii) holds. 
This completes the proof of Step \ref{s2 unbounded p=3}. 
\end{proof}

By extending $D$ to the fractional field $D: k'(x_0, x_1, x_2) \to k'(x_0, x_1, x_2)$ and restricting it to the
degree-$0$ and degree-$1$ parts, we obtain a $k$-derivation
\[
D:\MO_Y\to\MO_Y(1).
\]
Let $f:Y\to Y^D=:X$ be the quotient by $D$.

\begin{step}\label{s3 unbounded p=3}
There are finite surjective morphisms
\[
Y=\P^2_{k'}
\xrightarrow{f}
X
\xrightarrow{g}
\P^2_k
\]
of regular projective surfaces such that the composition is the base
change obtained by $k'/k$ and $f:Y\to X$ is an $\alpha_{-A}$-torsor for
$A:=g^*\MO_{\P^2_k}(1)$.
\end{step}

\begin{proof}[Proof of Step \ref{s3 unbounded p=3}]
The derivation $D$ is nowhere vanishing, since
$D(\alpha_0)=x_0$, $D(\alpha_1)=x_1$, $D(\alpha_2)=x_2$,
and $x_0,x_1,x_2$ have no common zero.
Since $D(k[x_0,x_1,x_2])=0$, the base-change morphism
$Y=\P^2_{k'} \to\P^2_k$ factors as
\[
Y=\P^2_{k'}
\xrightarrow{f}
X
\xrightarrow{g}
\P^2_k.
\]
Then $\MO_Y(1) = f^*g^*\MO_{\P^2_k}(1) =f^*\MO_{X}(A)$.
By Proposition \ref{p derivation to torsor}, $f$ is an
$\alpha_{-A}$-torsor. In particular, $f$ is faithfully flat, 
and hence $X$ is a regular projective surface.
This completes the proof of Step \ref{s3 unbounded p=3}. 
\end{proof}

\begin{step}\label{s4 unbounded p=3}
$-K_X$ is ample, $K_X^2 = 3^{n+2}$, and $-\chi(X, \MO_X) = 3^{n+1}$. 
\end{step}

\begin{proof}[Proof of Step \ref{s4 unbounded p=3}]
Assume that $H^0(X,\MO_{X})=k$. Since $f$ is an
$\alpha_{-A}$-torsor, we have
\[
\omega_Y
\simeq
f^*
\left(
\omega_{X}\otimes\MO_{X}(-2A)
\right).
\]
As $\omega_Y\simeq\MO_Y(-3)$ and
$f^*\MO_{X}(A)\simeq\MO_Y(1)$,
we obtain $f^*\omega_{X}\simeq\MO_Y(-1)$.
Thus $-K_{X}\equiv A$. In particular, $-K_X$ is ample. 
We have
\[
\deg g=
\frac{\deg (\P^2_{k'} \to \P^2_k)}{\deg f}
=
\frac{[k':k]}{\deg f}
=
\frac{3^{n+3}}{3}
=
3^{n+2}.
\]
Hence
$K_{X}^2
=
A^2
=
\deg g
=
3^{n+2}$.
As $f:Y\to X$ is an $\alpha_{-A}$-torsor, 
Proposition \ref{p torsor} implies 
\[
3^{n+3}=[k':k] = \chi(Y, \MO_Y) =
\chi(X,\MO_{X})
+
\chi(X,\MO_{X}(A))
+
\chi(X,\MO_{X}(2A)) 
\]
\[
= 3\chi(X, \MO_{X}) + 4K_X^2
= 3\chi(X, \MO_{X}) + 4 \cdot 3^{n+2}. 
\]
Therefore, $-\chi(X, \MO_X) = 3^{n+1}$. 
This completes the proof of Step \ref{s4 unbounded p=3}. 
\end{proof}

\begin{step}\label{s5 unbounded p=3}
For every $0\leq i\leq 2$, consider the $k$-derivation
$P_i:k'\to k'$ defined by
\[
D(\theta)
=
\sum_{i=0}^2P_i(\theta)x_i
\in
k'x_0\oplus k'x_1\oplus k'x_2.
\]
Set
$C
:=
[P_0,P_1]
=
P_0P_1-P_1P_0$, which is a $k$-derivation.
Then
$C^{3^{r-1}}|_{k_r}
=
\partial_{\beta_r}$
and $C^{3^r}|_{k_r}=0$
for every $1\leq r\leq n$.
\end{step}

\begin{proof}[Proof of Step \ref{s5 unbounded p=3}]
Recall that 
\[
D(\alpha_i)= x_i, \qquad 
D(\beta_r) = F_rx_1, \qquad 
F_{r+1}
=
-\beta_r^2F_r+F_r^3\alpha_1^2.
\]
Hence 
$P_0
=
\partial_{\alpha_0}$,
$P_1
=
\partial_{\alpha_1}
+
\sum_{r=1}^nF_r\partial_{\beta_r}$,
$P_2
=
\partial_{\alpha_2}$.
In particular, $P_0(\beta_r)=0$
for every $1\leq r\leq n$.
We have
\[
C(\beta_r) = (P_0P_1 - P_1P_0)(\beta_r) = P_0(P_1(\beta_r)) - P_1(P_0(\beta_r)) = P_0(F_r) 
\]
for every $1\leq r\leq n$. For $1\leq r<n$, we obtain 
\[
\begin{aligned}
C(\beta_{r+1})
&=
P_0(F_{r+1})\\
&=
P_0\left(
-\beta_r^2F_r+F_r^3\alpha_1^2
\right)\\
&=
-\beta_r^2 P_0(F_r) + F_r^3P_0(\alpha_1^2)\\
&=
-\beta_r^2C(\beta_r).
\end{aligned}
\]


We prove the assertion by induction on $r$. 
The base case $r=1$ of this induction holds, because 
the equalities 
$P_0|_{k_1}
=
\partial_{\alpha_0}$
and
$P_1|_{k_1}
=
\partial_{\alpha_1}
+
\alpha_0\partial_{\beta_1}$ give
\[
C|_{k_1} = [P_0|_{k_1}, P_1|_{k_1}] = 
[\partial_{\alpha_0}, \partial_{\alpha_1}
+
\alpha_0\partial_{\beta_1}]=[\partial_{\alpha_0}, \alpha_0\partial_{\beta_1}]
=\partial_{\beta_1}, 
\]
which implies $C^3|_{k_1} = \partial_{\beta_1}^3 =0$. 

Fix $1\leq r<n$ and assume that the assertion holds for $r$. Then
\[
C^{3^{r-1}}|_{k_r}
=
\partial_{\beta_r},
\qquad
C^{3^r}|_{k_r}=0,
\qquad
C^{3^{r-1}}(\beta_r)=1,
\qquad
C^{3^{r-1}+1}(\beta_r)=0.
\]
It suffices to show $C^{3^r}(\beta_{r+1})=1$, 
because 
$C^{3^r}|_{k_r}=0$ and $C^{3^r}(\beta_{r+1})=1$ give  $C^{3^r}|_{k_{r+1}}=\partial_{\beta_{r+1}}$, which implies  $C^{3^{r+1}}|_{k_{r+1}}=0$. 
It holds that 
\[
\begin{aligned}
C^{3^r}(\beta_{r+1})
&=
C^{3^r-1}
\left(
C(\beta_{r+1})
\right)\\
&=
C^{3^r-1}
\left(
-\beta_r^2C(\beta_r)
\right)\\
&\overset{(\star)}{=}
-
\sum_{\substack{a+b+c=3^r-1\\a,b,c\geq0}}
\binom{3^r-1}{a,b,c}
C^a(\beta_r)
C^b(\beta_r)
C^{c+1}(\beta_r)\\
&\overset{(\star\star)}{=}
-(-1)
C^{3^{r-1}}(\beta_r)
C^{3^{r-1}}(\beta_r)
C^{3^{r-1}}(\beta_r)\\
&=
1.
\end{aligned}
\]
Here $(\star)$ holds by the Leibniz rule and $(\star\star)$ follows from 
$C^d(\beta_r) =0$ for $d>3^{r-1}$ and the congruence 
\[
\binom{3^r-1}
{3^{r-1},3^{r-1},3^{r-1}-1}
\equiv \binom{2}
{1,1,0}
\times  
\binom{2}{0,0,2}^{r-1} =2 \equiv 
-1
\pmod 3, 
\]
which is ensured by the multinomial Lucas theorem. 
This completes the proof of Step \ref{s5 unbounded p=3}. 
\qedhere




\end{proof}

\begin{step}\label{s6 unbounded p=3}
It holds that $H^0(X,\MO_{X})=k$.
\end{step}

\begin{proof}[Proof of Step \ref{s6 unbounded p=3}]
Take $\theta\in k'=k_n$ satisfying $D(\theta)=0$. 
It suffices to show $\theta \in k$. By
\[
0
=
D(\theta)
=
\sum_{i=0}^2P_i(\theta)x_i
\in
k'x_0\oplus k'x_1\oplus k'x_2,
\]
we obtain $P_i(\theta)=0$
for every $0\leq i\leq 2$. In particular,
$C(\theta)
=
[P_0,P_1](\theta)
=
0$.
Hence
$\partial_{\beta_n}(\theta)
=
C^{3^{n-1}}(\theta)
=
0$, which implies $\theta\in k_{n-1}$. 
Next, we obtain $\partial_{\beta_{n-1}}(\theta)
=
(C^{3^{n-2}}|_{k_{n-1}})(\theta)
=C^{3^{n-2}}(\theta)=
0$, which implies $\theta\in k_{n-2}$. 
Repeating this procedure, we obtain
$\theta\in k_0
=
k(\alpha_0,\alpha_1,\alpha_2)$.
On $k_0$, we have
\[
D(\theta)
=
\partial_{\alpha_0}(\theta)x_0
+
\partial_{\alpha_1}(\theta)x_1
+
\partial_{\alpha_2}(\theta)x_2.
\]
Since $x_0,x_1,x_2$ are linearly independent over $k_0$, we obtain
$\partial_{\alpha_0}(\theta)
=
\partial_{\alpha_1}(\theta)
=
\partial_{\alpha_2}(\theta)
=
0$.
Therefore, $\theta\in k$. Thus
$H^0(X,\MO_{X})=k$.
This completes the proof of Step \ref{s6 unbounded p=3}. 
\end{proof}
Step \ref{s4 unbounded p=3} and Step \ref{s6 unbounded p=3}
 complete the proof of Theorem \ref{t unbounded p=3}. 
\end{proof}

\begin{cor}\label{c unbounded p=3}
Let $k$ be a field of characteristic three satisfying $[k:k^3] =\infty$. 
Then, for every integer $L>0$, there exists a regular del Pezzo surface $X$ over $k$ such that $K_X^2 \geq L$ and $h^1(X, \MO_X) \geq L$. 
\end{cor}

\begin{proof}
The assertion follows from Theorem \ref{t unbounded p=3}.
\end{proof}

\subsection{Case $p=2$}

\begin{prop}\label{p unified p=2 construction}
Let $k$ be a field of characteristic two. 
Take integers $m,n\geq1$ and an intermediate field 
\[
k
\subset
k'
:=
k(\alpha_0,\ldots,\alpha_{m-1},\beta_1,\ldots,\beta_n)
\subset
k^{1/2}
\]
satisfying $[k':k]=2^{m+n}$. 
Let $V$ be a regular projective surface over $k$ 
such that $H^0(V,\MO_V)=k$ and $Y:=V_{k'}$ is regular. 
Take a Cartier divisor $L_V$ on $V$ and linearly independent sections $x_0,\ldots,x_{m-1}\in H^0(V,\MO_V(L_V))$ having no common zero. 
Assume that 
$H_V:=-K_V-L_V$ is ample and 
$h^0(V,\MO_V(2L_V)) <\frac{m(m+1)}{2}$. 
Then there exists a regular del Pezzo surface $X$ over $k$ such that 
\begin{enumerate}
\item $K_X^2=2^{m+n-1}H_V^2$, 
\item $-\chi(X, \MO_X) = 2^{m+n-3}L_V\cdot (-K_V)
-2^{m+n-1}\chi(V,\MO_V)$. 
\item 
There are finite surjective morphisms
\[
Y=V_{k'}\xrightarrow{f}X\xrightarrow{g}V
\]
such that the composition is the base-change morphism obtained by $k'/k$,
$f$ is a non-trivial $\alpha_{-M}$-torsor for $M:=g^*L_V$, and 
$(X\times_k\overline{k})_{\red}^N \simeq (V\times_k\overline{k})_{\red}^N$. 
\end{enumerate}
\end{prop}

\begin{proof}
For $0\leq r\leq n$, set
$k_r:=k(\alpha_0,\ldots,\alpha_{m-1},\beta_1,\ldots,\beta_r)$. 
For example,
\[
k_0=k(\alpha_0,\ldots,\alpha_{m-1})
\qquad\text{and}\qquad
k_n=k(\alpha_0,\ldots,\alpha_{m-1},\beta_1,\ldots,\beta_n)=k'.
\]
Let $L$ and $H$ be the pullbacks of $L_V$ and $H_V$ to $Y$, respectively. 
We use the same notation $x_i$ for the pullbacks of $x_i$ to $Y$.

\setcounter{step}{0}

\begin{step}\label{s0 unified p=2}
After replacing $x_0, ..., x_{m-1} \in H^0(V,\MO_V(L_V))$ by 
a $k$-linear basis of $kx_0 \oplus \cdots \oplus kx_{m-1}$ if necessary, 
there exist $d_i,c_{ij}\in k$ 
satisfying 
\[
x_0^2
+
\sum_{i=1}^{m-1}d_ix_0x_i
+
\sum_{1\leq i\leq j\leq m-1}c_{ij}x_ix_j
=
0 \qquad \text{in} \qquad H^0(V, \MO_V(2L_V)).  
\]
\end{step}
\begin{proof}[Proof of Step \ref{s0 unified p=2}]
Set
$W:= k x_0 \oplus \cdots \oplus k x_{m-1} 
\subset H^0(V,\MO_V(L_V))$.
Consider the multiplication map
\[
\mu:\Sym^2W\longrightarrow H^0(V,\MO_V(2L_V)).
\]
Since $\dim_k\Sym^2W=m(m+1)/2>h^0(V,\MO_V(2L_V))$, 
we have $\Ker\mu\neq0$. Thus there exist elements
$a_{ij}\in k$ such that some $a_{ij}$ is nonzero and 
\[
\sum_{0\leq i\leq j\leq m-1}a_{ij}x_ix_j=0
\qquad\text{in}\qquad
H^0(V,\MO_V(2L_V)).
\]
If $a_{ii} \neq 0$ for some $i$, 
then we may assume that $a_{00} \neq 0$, and hence we are done by 
taking the multiplication with $a_{00}^{-1}$. 
Otherwise, there exist $i<j$ such that 
$a_{ij}\neq0$. After reordering, we may assume that $(i,j)=(0,1)$.
Then the problem is reduced to the case when $a_{00} \neq 0$ 
by replacing the basis $x_0,x_1,x_2,\ldots,x_{m-1}$ by
$x_1,x_0+x_1,x_2,\ldots,x_{m-1}$. 
This completes the proof of Step \ref{s0 unified p=2}.
\qedhere


\end{proof}
\begin{step}\label{s1 unified p=2}
There exists an $\MO_V$-linear derivation
\[
D:\MO_Y
\to
\MO_Y(L)
\]
satisfying the following conditions:
\begin{enumerate}
\item $D(\alpha_i)=x_i$ for every $0\leq i\leq m-1$.
\item $D(k_r)\subset\bigoplus_{i=0}^{m-1}k_rx_i$ for every $0\leq r\leq n$.
\item It holds that 
\[
D(\beta_1)
=
\alpha_0x_0
+
\sum_{i=1}^{m-1}d_i\alpha_0x_i
+
\sum_{1\leq i\leq j\leq m-1}c_{ij}\alpha_ix_j.
\]
\item If $1\leq r<n$ and
$D(\beta_r)=\sum_{i=0}^{m-1}F_{i,r}x_i$ for $F_{i,r}\in k_r$, then
\[
D(\beta_{r+1})
=
\beta_rD(\beta_r)
+
\sum_{i=0}^{m-1}F_{i,r}^2\alpha_ix_i.
\]
\end{enumerate}
\end{step}

\begin{proof}[Proof of Step \ref{s1 unified p=2}]
By $\MO_Y = \MO_{V_{k'}} =\MO_V\otimes_k k'$, it is enough to find a $k$-derivation
$D:k'\to\bigoplus_{i=0}^{m-1}k'x_i$ satisfying (1)-(4).
There exists a unique $k$-derivation
$D_0:k_0\to\bigoplus_{i=0}^{m-1}k_0x_i$ satisfying (1).
Then there is a unique $k$-derivation
$D_1:k_1\to\bigoplus_{i=0}^{m-1}k_1x_i$ such that
$D_1|_{k_0}=D_0$ and (3) holds.
Similarly, we can inductively extend
$D_r:k_r\to\bigoplus_{i=0}^{m-1}k_rx_i$ to
$D_{r+1}:k_{r+1}\to\bigoplus_{i=0}^{m-1}k_{r+1}x_i$
by using (4).
Then the resulting $k$-derivation
$D_n:k_n\to\bigoplus_{i=0}^{m-1}k_nx_i$ satisfies (1)-(4).
This completes the proof of Step \ref{s1 unified p=2}.
\end{proof}

\begin{step}\label{s2 unified p=2}
It holds that $D^2=0$.
\end{step}

\begin{proof}[Proof of Step \ref{s2 unified p=2}]
For every $0\leq i\leq m-1$, we have
$D^2(\alpha_i)\overset{{\rm (1)}}{=}D(x_i)=0$.
As $D^2$ is a $k$-derivation, it suffices to show that
$D^2(\beta_r)=0$ for every $1\leq r\leq n$.
We prove this by induction on $r$.
The base case $r=1$ follows from
\[
\begin{aligned}
D^2(\beta_1)
&\overset{{\rm (3)}}{=}
D\left(
\alpha_0x_0
+
\sum_{i=1}^{m-1}d_i\alpha_0x_i
+
\sum_{1\leq i\leq j\leq m-1}c_{ij}\alpha_ix_j
\right)\\
&\overset{{\rm (1)}}{=}
x_0^2
+
\sum_{i=1}^{m-1}d_ix_0x_i
+
\sum_{1\leq i\leq j\leq m-1}c_{ij}x_ix_j\\
&=
0.
\end{aligned}
\]
Fix $1\leq r<n$ and assume that $D^2(\beta_r)=0$. Then
\[
\begin{aligned}
D^2(\beta_{r+1})
&\overset{{\rm (4)}}{=}
D\left(
\beta_rD(\beta_r)
+
\sum_{i=0}^{m-1}F_{i,r}^2\alpha_ix_i
\right)\\
&=
D(\beta_r)^2
+
\beta_rD^2(\beta_r)
+
\sum_{i=0}^{m-1}F_{i,r}^2x_i^2\\
&=
\left(
\sum_{i=0}^{m-1}F_{i,r}x_i
\right)^2
+0+
\sum_{i=0}^{m-1}F_{i,r}^2x_i^2\\
&=
0.
\end{aligned}
\]
Therefore, $D^2(\beta_r)=0$ for every $1\leq r\leq n$.
Thus $D^2=0$.
This completes the proof of Step \ref{s2 unified p=2}.
\end{proof}

Let $f:Y\to Y^D=:X$ be the quotient by $D$.

\begin{step}\label{s3 unified p=2}
There are finite surjective morphisms
\[
Y=V_{k'}\xrightarrow{f}X\xrightarrow{g}V
\]
of regular projective surfaces such that 
$(X\times_k\overline{k})_{\red}^N
\xrightarrow{\simeq} 
(V\times_k\overline{k})_{\red}^N$, the composition is
the base-change morphism obtained by $k'/k$, and
$f:Y\to X$ is a non-trivial $\alpha_{-M}$-torsor for $M:=g^*L_V$.
\end{step}

\begin{proof}[Proof of Step \ref{s3 unified p=2}]
The derivation $D:\MO_Y
\to
\MO_Y(L)$ is nowhere vanishing, because
\[
D(\alpha_0)=x_0,\qquad
D(\alpha_1)=x_1,\qquad\ldots,\qquad
D(\alpha_{m-1})=x_{m-1}
\]
and $x_0,\ldots,x_{m-1}$ have no common zero.
Since $D(\MO_V)=0$, the base-change morphism $Y=V_{k'}\to V$
factors as
\[
Y=V_{k'}\xrightarrow{f}X\xrightarrow{g}V.
\]
Then $\MO_Y(L)=f^*g^*\MO_V(L_V)=f^*\MO_X(M)$.
By Proposition \ref{p derivation to torsor}, $f$ is an
$\alpha_{-M}$-torsor.
In particular, $f$ is faithfully flat, and hence $X$ is a regular
projective surface.

Since $k'/k$ is purely inseparable, we have
\[
(Y\times_k\overline{k})_{\red}=( (V \times_k k') \times_k \ol{k})_{\red} 
\xrightarrow{\iota, \simeq} 
(V\times_k\overline{k})_{\red}.
\]
Then the induced morphisms 
\[
\iota^N: (V\times_k\overline{k})^N_{\red}
\to
(X\times_k\overline{k})^N_{\red}
\to
(V\times_k\overline{k})^N_{\red}
\]
are finite surjective morphisms of projective normal varieties over $\ol{k}$ 
whose composition $\iota^N$ is an isomorphism. 
Thus $(X\times_k\overline{k})^N_{\red}
\to
(V\times_k\overline{k})^N_{\red}$ is also an isomorphism. 
This completes the proof of Step \ref{s3 unified p=2}.
\end{proof}

\begin{step}\label{s4 unified p=2}
$-K_X$ is ample, 
$K_X^2
=
2^{m+n-1}H_V^2$, and 
\[
-\chi(X, \MO_X) = 2^{m+n-3}(-L_V\cdot K_V)-
2^{m+n-1}\chi(V,\MO_V).
\]
\end{step}

\begin{proof}[Proof of Step \ref{s4 unified p=2}]
Since $f$ is an $\alpha_{-M}$-torsor, we have
\[
\omega_Y
\simeq
f^*
\left(
\omega_X\otimes\MO_X(-M)
\right).
\]
As $\omega_Y\simeq(g\circ f)^*\omega_V$ and
$f^*\MO_X(M)\simeq\MO_Y(L)\simeq(g\circ f)^*\MO_V(L_V)$, we obtain
\[
f^*K_X
\sim
(g\circ f)^*(K_V+L_V)
=
-(g\circ f)^*H_V.
\]
Thus $-K_X\equiv g^*H_V$, which is ample. 
We have
\[
\deg g
=
\frac{\deg(V_{k'}\to V)}{\deg f}
=
\frac{[k':k]}{2}
=
2^{m+n-1}.
\]
Hence
$K_X^2
=
(g^*H_V)^2
=
(\deg g)H_V^2
=
2^{m+n-1}H_V^2.$ 

As $f:Y\to X$ is an $\alpha_{-M}$-torsor, Proposition \ref{p torsor}(3) gives
\[
\chi(Y,\MO_Y)
=
\chi(X,\MO_X)
+
\chi(X,\MO_X(M)).
\]
By $-K_X\equiv g^*H_V$ and the Riemann--Roch theorem, we obtain
\[
\begin{aligned}
\chi(X,\MO_X(M))
&=
\chi(X,\MO_X)
+
\frac{1}{2}M\cdot(M-K_X)\\
&=
\chi(X,\MO_X)
+
\frac{\deg g}{2}L_V\cdot(L_V+H_V)\\
&=
\chi(X,\MO_X)
+
2^{m+n-2}L_V\cdot (-K_V).
\end{aligned}
\]
On the other hand,
$\chi(Y,\MO_Y)=[k':k]\chi(V,\MO_V)=2^{m+n}\chi(V,\MO_V)$.
To summarise, 
\[
-\chi(X,\MO_X)
=
2^{m+n-3}L_V\cdot (-K_V)
-2^{m+n-1}\chi(V,\MO_V).
\]
This completes the proof of Step \ref{s4 unified p=2}.
\end{proof}

\begin{step}\label{s5 unified p=2}
For every $0\leq i\leq m-1$, consider the $k$-derivation
\[
P_i:k'\to k'
\qquad\text{defined by}\qquad
D(\theta)
=
\sum_{i=0}^{m-1}P_i(\theta)x_i
\in
\bigoplus_{i=0}^{m-1}k'x_i.
\]
Set $P:=P_0$. Then
$P^{2^r}|_{k_r}=\partial_{\beta_r}$ and
$P^{2^{r+1}}|_{k_r}=0$ for every $1\leq r\leq n$.
\end{step}

\begin{proof}[Proof of Step \ref{s5 unified p=2}]
Recall that
\[
D(\alpha_i)=x_i,\qquad
D(\beta_r)=\sum_{i=0}^{m-1}F_{i,r}x_i,\qquad
D(\beta_{r+1})
=
\beta_rD(\beta_r)
+
\sum_{i=0}^{m-1}F_{i,r}^2\alpha_ix_i, 
\]
\[
D(\beta_1)
=
\alpha_0x_0
+
\sum_{i=1}^{m-1}d_i\alpha_0x_i
+
\sum_{1\leq i\leq j\leq m-1}c_{ij}\alpha_ix_j.
\]
Hence $F_{0, r} = P(\beta_r)$, $F_{0, 1} = P(\beta_1) = \alpha_0$, and 
\[
\sum_{i=0}^{m-1}F_{i,r+1}x_i
=D(\beta_{r+1})
=
\beta_rD(\beta_r)
+
\sum_{i=0}^{m-1}F_{i,r}^2\alpha_ix_i=
\sum_{i=0}^{m-1}
\left(
\beta_rF_{i,r}+F_{i,r}^2\alpha_i
\right)x_i.
\]
The equality of the coefficients of $x_i$ gives
$F_{i,r+1}=\beta_rF_{i,r}+F_{i,r}^2\alpha_i$.
For $f_r:=F_{0,r}=P(\beta_r)$, 
we get $f_1 = F_{0, 1} = \alpha_0$ and 
$f_{r+1}
=
\beta_rf_r+f_r^2\alpha_0$ for every $1\leq r<n$.
Hence
\[
P
=
\partial_{\alpha_0}
+
\sum_{r=1}^nf_r\partial_{\beta_r}.
\]

We prove the assertion by induction on $r$.
The base case $r=1$ holds, because
$P|_{k_1}
=
\partial_{\alpha_0}
+
\alpha_0\partial_{\beta_1}$ gives
$P^2|_{k_1}=(\partial_{\alpha_0}
+
\alpha_0\partial_{\beta_1})(\partial_{\alpha_0}
+
\alpha_0\partial_{\beta_1})=\partial_{\beta_1}$ and
$P^4|_{k_1}=\partial_{\beta_1}^2 = 0$.

Fix $1\leq r<n$ and assume that the assertion holds for $r$. Then
\[
P^{2^r}|_{k_r}
=
\partial_{\beta_r},
\qquad
P^{2^{r+1}}|_{k_r}=0,
\qquad
P^{2^r}(\beta_r)=1,
\qquad
P^{2^r+1}(\beta_r)=0.
\]
It suffices to show $P^{2^{r+1}}(\beta_{r+1})=1$, because
$P^{2^{r+1}}|_{k_r}=0$ and
$P^{2^{r+1}}(\beta_{r+1})=1$ give
$P^{2^{r+1}}|_{k_{r+1}}=\partial_{\beta_{r+1}}$, which implies
$P^{2^{r+2}}|_{k_{r+1}}=0$.
We have $P(f_r^2)=0$ and $P^2(\alpha_0) =P(1)=0$. Hence
$P^{2^{r+1}-1}(f_r^2\alpha_0)=0$. It follows that
\[
\begin{aligned}
P^{2^{r+1}}(\beta_{r+1})
&=
P^{2^{r+1}-1}
\left(
P(\beta_{r+1})
\right)\\
&=
P^{2^{r+1}-1}(f_{r+1})\\
&=
P^{2^{r+1}-1}
\left(
\beta_rf_r+f_r^2\alpha_0
\right)\\
&=
P^{2^{r+1}-1}(\beta_rf_r)\\
&\overset{(\star)}{=}
\sum_{i=0}^{2^{r+1}-1}
P^i(\beta_r)
P^{2^{r+1}-i}(\beta_r)\\
&\overset{(\star\star)}{=}
P^{2^r}(\beta_r)
P^{2^r}(\beta_r)\\
&=
1.
\end{aligned}
\]
Here $(\star)$ follows from the Leibniz rule and the congruence
$\binom{2^{r+1}-1}{i}\equiv1\pmod2$ for every
$0\leq i\leq2^{r+1}-1$, which follows from
$(1+T)^{2^{r+1}}=1+T^{2^{r+1}}
=
(1+T)(1+T+\cdots+T^{2^{r+1}-1})$
in $\F_2[T]$.
The equality $(\star\star)$ holds by
$P^d(\beta_r)=0$ for every $d>2^r$.
This completes the proof of Step \ref{s5 unified p=2}.
\qedhere
\end{proof}

\begin{step}\label{s6 unified p=2}
It holds that $H^0(X,\MO_X)=k$.
\end{step}

\begin{proof}[Proof of Step \ref{s6 unified p=2}]
Take $\theta\in k'=k_n$ satisfying $D(\theta)=0$. It suffices to show
$\theta\in k$.
By
\[
0
=
D(\theta)
=
\sum_{i=0}^{m-1}P_i(\theta)x_i
\in
\bigoplus_{i=0}^{m-1}k'x_i,
\]
we obtain $P_i(\theta)=0$ for every $0\leq i\leq m-1$.
In particular, $P(\theta) =P_0(\theta)=0$.
Hence
$\partial_{\beta_n}(\theta)
=
(P^{2^n}|_{k_n})(\theta)
=
0$,
and therefore $\theta\in k_{n-1}$.
We then obtain
$\partial_{\beta_{n-1}}(\theta)
=
(P^{2^{n-1}}|_{k_{n-1}})(\theta)
=
0$.
Hence $\theta\in k_{n-2}$.
Repeating this procedure, we obtain $\theta\in k_0$.

It holds that
\[
0
=
D(\theta)
=
\sum_{i=0}^{m-1}
\partial_{\alpha_i}(\theta)x_i
\in
\bigoplus_{i=0}^{m-1}k_0x_i,
\]
which implies $\partial_{\alpha_i}(\theta)=0$ for every
$0\leq i\leq m-1$.
This, together with
$k_0=k(\alpha_0,\ldots,\alpha_{m-1})$ and
$[k_0:k]=2^m$, implies $\theta\in k$.
This completes the proof of Step \ref{s6 unified p=2}.
\end{proof}
Steps \ref{s3 unified p=2}, \ref{s4 unified p=2}, and
\ref{s6 unified p=2} complete the proof of Proposition
\ref{p unified p=2 construction}.
\end{proof}

\begin{thm}\label{t unbounded p=2}
Take an integer $n\geq1$ and a field $k$ of characteristic two satisfying
$[k:k^2]\geq2^{n+6}$.
Then there exists a regular del Pezzo surface $X$ over $k$ such that
$K_X^2=2^{n+5}$ and $-\chi(X, \MO_X) = 2^{n+4}$. 
\end{thm}

\begin{proof}
By $[k^{1/2}:k] =[k:k^2]\geq2^{n+6}$, we can find elements
$\alpha_0,\ldots,\alpha_5,\beta_1,\ldots,\beta_n\in k^{1/2}$ such that the 
intermediate field 
\[
k
\subset
k'
:=
k(\alpha_0,\ldots,\alpha_5,\beta_1,\ldots,\beta_n)
\subset
k^{1/2}
\]
satisfies 
$[k':k]=2^{n+6}$. 

Set $V:=\P^2_k = \Proj k [z_0, z_1, z_2]$, $L_V:=\MO_{\P^2_k}(2)$, and 
\[
x_0:=z_0z_1,\quad
x_1:=z_0^2,\quad
x_2:=z_1^2,\quad
x_3:=z_2^2,\quad
x_4:=z_0z_2,\quad
x_5:=z_1z_2 \in H^0(V, L_V).
\]
Then $x_0,\ldots,x_5$ form a $k$-linear basis of
$H^0(V,\MO_V(L_V))$ and have no common zero. 
Moreover,
$H_V:=-K_V-L_V=\MO_{\P^2_k}(1)$,
$H_V^2=1$, $L_V\cdot (-K_V)=6$, and $\chi(V,\MO_V)=1$. 
Set $m:=6$. 
Since $Y:=V_{k'}\simeq\P^2_{k'}$ is regular and 
$m(m+1)/2 = 21 > 15 = h^0(\P^2, \MO_{\P^2}(4))=h^0(V, \MO_V(2L_V))$,
Proposition \ref{p unified p=2 construction} is applicable. Therefore, $K_X^2=2^{m+n-1} = 2^{n+5}$ and
\[
-\chi(X, \MO_X) = 2^{m+n-3}L_V\cdot (-K_V)
-2^{m+n-1}\chi(V,\MO_V) = 2^{n+3} \cdot 6 - 2^{n+5} \cdot 1=2^{n+4}, 
\]
as required. 
\end{proof}

\begin{cor}\label{c unbounded p=2}
Let $k$ be a field of characteristic two satisfying $[k:k^2] =\infty$. 
Then, for every integer $L>0$, there exists a regular del Pezzo surface $X$ over $k$ such that $K_X^2 \geq L$ and $h^1(X, \MO_X) \geq L$. 
\end{cor}

\begin{proof}
The assertion follows from Theorem \ref{t unbounded p=2}.
\end{proof}

\begin{thm}\label{t unbounded pathological}
For every integer $n\geq1$, there exist a field $k$ of characteristic two
and a regular del Pezzo surface $X$ over $k$ such that
$K_X^2=2^{n+8}$ and
$-\chi(X, \MO_X) =2^{n+6}$.
\end{thm}

\begin{proof}
By Example \ref{e index2 K^2=16}, 
we can find a field $\bF$ of characteristic two,
a regular del Pezzo surface $V_0$ over $\bF$, and an ample Cartier divisor $L_{V_0}$ on $V_0$ such that
$H^0(V_0,\MO_{V_0})=\bF$, $-K_{V_0}\sim2L_{V_0}$, $L_{V_0}^2=4$,
$h^0(V_0,L_{V_0})=7$, $h^0(V_0,2L_{V_0})=17$, and $|L_{V_0}|$ is base point free. 
Take an $\bF$-linear basis
$x_0,\ldots,x_6$ of $H^0(V_0,\MO_{V_0}(L_{V_0}))$.
Since $|L_{V_0}|$ is base point free, the elements
$x_0,\ldots,x_6$ have no common zero. 

Take algebraically independent elements
$a_0,\ldots,a_6,b_1,\ldots,b_n$ over $\bF$ and set
\[
k
:=
\bF(a_0,\ldots,a_6,b_1,\ldots,b_n)
\subset
k(\alpha_0,\ldots,\alpha_6,\beta_1,\ldots,\beta_n)
=:k',
\]
where $\alpha_i^2=a_i$ and $\beta_r^2=b_r$.
Then $[k':k]=2^{n+7}$.
Set $V:=V_0\times_{\bF}k$ and $Y:=V_0\times_{\bF}k'=V_{k'}$, which are regular projective surfaces. 
For the pullback $L_V$ of $L_{V_0}$ on $V$, 
we set $H_V := -K_V -L_V \sim L_V$. 
Then $H_V$ is ample, $H_V^2 = L_V^2 =4$, 
and $\chi(V, \MO_V) = 1$. 
By 
\[
\frac{7\cdot8}{2}
>
17
=
h^0(V,\MO_{V}(2L_V)), 
\]
Proposition \ref{p unified p=2 construction} is applicable for $m:=7$. 
Hence there exists a regular del Pezzo surface $X$ over $k$ 
such that $K_X^2
=
2^{n+6}\cdot H_V^2 =
2^{n+8}$ 
and
\[
-\chi(X, \MO_X) = 2^{m+n-3}L_V\cdot (-K_V)
-2^{m+n-1}\chi(V,\MO_V) = 2^{n+4} \cdot 8 -2^{n+6} \cdot 1 
=2^{n+6}. 
\]
\end{proof}

\begin{thm}\label{t unbounded h1=1}
For every integer $n\geq1$ and every 
$Z\in\{\P^1\times\P^1,\P(1,1,2),\P^2\}$, 
there exist a field $k$ of characteristic two
and a regular del Pezzo surface $X$ over $k$ such that
$K_X^2=2^{n+4}$, $\chi(X, \MO_X)=0$, and
$(X\times_k\overline{k})_{\red}^N\simeq Z$.
\end{thm}

\begin{proof}
Let $\bF$ be a field of characteristic two  and take
$u,v,w\in\bF$ such that 
$[\bF^2(u,v,w):\bF^2]=8$.
In $\P^3_{\bF} = \Proj \bF[z_0, z_1, z_2, z_3]$, set
\[
\begin{aligned}
W_1&:=\{z_0z_1+z_2z_3=0\},\\
W_2&:=\{z_0z_1+z_2^2+uz_3^2=0\},\\
W_3&:=\{z_0^2+uz_1^2+vz_2^2+wz_3^2=0\}.
\end{aligned}
\]
Then $W_1,W_2,W_3$ are regular quadric surfaces. Moreover,
\[
(W_1\times_{\bF}\overline{\bF})_{\red}^N
\simeq\P^1\times\P^1,\qquad
(W_2\times_{\bF}\overline{\bF})_{\red}^N
\simeq\P(1,1,2),\qquad 
(W_3\times_{\bF}\overline{\bF})_{\red}^N\simeq\P^2.
\]

Fix an integer $n\geq1$. Take algebraically independent elements
$a_0,\ldots,a_3,b_1,\ldots,b_n$ over $\bF$ and set 
\[
k:=\bF(a_0,\ldots,a_3,b_1,\ldots,b_n) \subset k(\alpha_0,\ldots,\alpha_3,\beta_1,\ldots,\beta_n)=:k', 
\]
where $\alpha_i^2=a_i$ and $\beta_r^2=b_r$. Then
$[k':k]=2^{n+4}$. 
For each $i\in\{1,2,3\}$, set
$V_i:=(W_i)_k$, $Y:=(V_i)_{k'}$, $L_{V_i}:=\MO_{V_i}(1)$, and $H_{V_i}:=-K_{V_i}-L_{V_i} \sim L_{V_i}$. Then $H_{V_i}$ is ample and
\[
H_{V_i}^2=2,\qquad
-L_{V_i}\cdot K_{V_i}=4,\qquad
\chi(V_i,\MO_{V_i})=1.
\]
The restrictions $x_0,\ldots,x_3 \in H^0(V_i,\MO_{V_i}(L_{V_i}))$ 
of $z_0,\ldots,z_3 \in H^0(\P^3_k, \MO_{\P^3_k}(1))$ are linearly independent and have no common zero. 
For  $m :=4$, we have 
$h^0(V_i,\MO_{V_i}(2L_{V_i}))=9<10=\frac{m\cdot(m+1)}{2}$.
Then Proposition \ref{p unified p=2 construction} is applicable. 
For every $i \in \{1, 2, 3\}$, we obtain a regular
del Pezzo surface $X_i$ over $k$ 
such that 
$(X_i\times_k\overline{k})_{\red}^N
\simeq
(V_i\times_k \overline{k})_{\red}^N \simeq (W_i \times_{\bF} \overline{k})_{\red}^N$, 
$K_{X_i}^2
=
2^{n+3}H_{V_i}^2
=
2^{n+4}$, 
and
\[
-\chi(X_i, \MO_{X_i})
=2^{n+1}L_{V_i}\cdot (-K_{V_i})
-2^{n+3}\chi(V_i,\MO_{V_i})=
4 \cdot 2^{n+1} -2^{n+3} =0. 
\]
\end{proof}

\subsection{Summary of examples with $H^1(X, \MO_X) \neq 0$}\label{ss irregular summary}

The following table summarises the four cases in Theorem \ref{t structure} and records the values of $K_X^2$ realised by the examples constructed in this paper.


\[
{\setlength{\arraycolsep}{3pt}
\begin{array}{c|c|c|c|c|c}
p
& \rho(X)
& \displaystyle-\frac{\chi(X,\MO_X)}{K_X^2}
& K_X^2(\text{existence})
& (X \times_k \overline{k})_{\red}^N
& \begin{array}{c}\text{Torsor description}\end{array}\\
\hline
2 & 1, 2 & 0
& 2^e\text{ (for } \forall e\geq 3)
& \begin{array}{c}
\P^2,\ \P^1\times\P^1,\\
\P(1,1,2)
\end{array}
& \begin{array}{c}
\exists\,\alpha_{\omega_X}\text{-torsor}\,\,Y\to X,\\
Y\text{:normal quadric}
\end{array}\\
\hline
2 & 1 & 1/4
& 2^e\text{ (for }  \forall e\geq 4)
& \P^2
& \begin{array}{c}
\exists\,\alpha_{\omega_X}\text{-torsor}\,\,Y\to X,\\
Y\text{:normal},\ K_Y^2=16[k_Y:k]
\end{array}\\
\hline
2 & 1 & 1/2
& 2^e\text{ (for }  \forall e\geq 2)
& \P^2
& \begin{array}{c}
\exists\,A\equiv-2K_X,\\
\exists\,\alpha_{-A}\text{-torsor}\,\,\P^2_{k_Y}\to X
\end{array}\\
\hline
3 & 1 & 1/3
& 3^e\text{ (for }  \forall e\geq 2)
& \P^2
& \begin{array}{c}
\exists\,\alpha_{\omega_X}\text{-torsor}\,\,\P^2_{k_Y}\to X
\end{array}\\
\hline
\end{array}}
\]
All the properties in the table except for the existence assertions follow from Theorem \ref{t structure}.
The existence is ensured as follows: 
\begin{itemize}
\item $-\chi(X, \MO_X)/K_X^2 = 0$: 
Example \ref{e rho=2 irregular} and Theorem \ref{t unbounded h1=1}.
\item $-\chi(X, \MO_X)/K_X^2 = 1/4$: 
Proposition \ref{p common pathological} and Theorem \ref{t unbounded pathological}
\item $-\chi(X, \MO_X)/K_X^2 = 1/2$: 
Example \ref{e H1(2KX)} and Theorem \ref{t unbounded p=2}
\item $-\chi(X, \MO_X)/K_X^2 = 1/3$: 
Example \ref{e p=3 irregular} and Theorem \ref{t unbounded p=3}. 
\end{itemize}

\begin{rem}
\begin{enumerate}
\item 
For each $i \in \{1, 2\}$, Maddock constructs 
a regular del Pezzo surface $X_i$ satisfying $K_{X_i}^2=i$ and $\chi(X_i, \MO_{X_i})=0$ \cite[Main Theorem(1)]{Mad16}. 
In particular, the only missing case for $\chi(X, \MO_X)=0$ is when $K_X^2=4$. 
\item 
If $X$ is geometrically integral, then we have $\chi(X, \MO_X)=0$ \cite[Theorem 3.7]{BT24}. 
In particular, we have $e \geq 1$ when $\chi(X, \MO_X) \neq 0$ (cf.\ \cite[Theorem 4.6]{Tan24}). 
\end{enumerate}
\end{rem}

\appendix

\section{Derivations and $\alpha_{\cL}$-torsors}

Given a derivation $D : \MO_Y \to \cL$ on $Y$, we get the quotient $f: Y \to Y^D =X$. 
The purpose of this appendix is to give a criterion for when $f$ is an $\alpha_{\cM}$-torsor (Proposition \ref{p derivation to torsor}). 
To this end, we start by establishing a couple of auxiliary results 
(Lemma \ref{l additive slice}, Lemma \ref{l unit slice}).

\begin{lem}\label{l additive slice}
Let $R$ be an $\mathbb F_p$-algebra and let
$D:R\to R$ be a derivation. Assume that $D^p=0$ and that there exists
$t\in R$ such that $D(t)=1$. Then
\[
R=\bigoplus_{i=0}^{p-1}R^Dt^i.
\]
In particular, $R$ is a free $R^D$-module of rank $p$.
\end{lem}

\begin{proof}
Let us prove that 
\[
R=\sum_{i=0}^{p-1} R^D t^i.
\]
Fix $r\in R$. We define $r_0, ..., r_{p-1} \in R$ and $a_0, ..., a_{p-1} \in R$ by the following recursive conditions: 
\begin{enumerate}
\item $r_{p-1} = r$. 
\item $a_n = \frac{1}{n!} D^n(r_n)$ if $0 \leq n \leq p-1$. 
\item $r_{n-1} =r_n -a_nt^n$ if $1 \leq n \leq p-1$. 
\end{enumerate}
Since (1)-(3) imply $r = \sum_{n=0}^{p-1} a_nt^n$, 
it is enough to show that $D(a_0) = D(a_1) = \cdots =D(a_{p-1})=0$. 
We have $D(a_{p-1}) = \frac{1}{(p-1)!} D^p(r_{p-1})=0$. 
If $1 \leq n \leq p-1$ and $D(a_n)=0$, then we get  
\[
(n-1)!D(a_{n-1})
\overset{{\rm (2)}}{=}
D^n(r_{n-1})
\overset{{\rm (3)}}{=}D^n(r_n -a_nt^n)
=
D^n(r_n)-a_nD^n(t^n)
=
D^n(r_n)-n!a_n
\overset{{\rm (2)}}{=}
0.
\]
This completes the proof of $R=\sum_{i=0}^{p-1} R^D t^i$.

It suffices to prove $\sum_{i=0}^{p-1} R^D t^i = \bigoplus_{i=0}^{p-1} R^D t^i$. 
Assume that 
$\sum_{i=0}^{p-1}b_it^i=0$ for some $b_0, ..., b_{p-1} \in R^D$. 
We have $0 = D^{p-1}(\sum_{i=0}^{p-1}b_it^i) = D^{p-1}(b_{p-1}t^{p-1}) = (p-1)!b_{p-1}$, which implies $b_{p-1}=0$. 
Then $0 = D^{p-2}(\sum_{i=0}^{p-2}b_it^i) = D^{p-2}(b_{p-2}t^{p-2}) = (p-2)!b_{p-2}$, which implies $b_{p-2}=0$. Repeating this procedure, we obtain $b_{p-1} = b_{p-2} = \cdots = b_0=0$, as required. 
\end{proof}

\begin{lem}\label{l unit slice}
Let $R$ be an $\mathbb F_p$-algebra and let
$D:R\to R$ be a derivation such that $D^p=0$.
Assume that there is $s\in R$ such that
$D(s)\in R^\times$. 
Then there exists $t\in R$ satisfying $D(t)=1.$ 
\end{lem}

\begin{proof}
Set $u:=D(s)\in R^\times$ and  
\[
E:=u^{-1}D : R \xrightarrow{D} R \xrightarrow{\times u^{-1}} R, 
\]
which is a derivation.
Then
$E(s)=1$ and $R^E=R^D$. 
By Hochschild's formula \cite[Theorem 25.5]{Mat86}, we get 
\[
E^p= (u^{-1}D)^p = u^{-p}D^p +(u^{-1}D)^{p-1}(u^{-1}) \cdot D = cE
\]
for $c := uE^{p-1}(u^{-1})\in R$. Then it holds that 
\[
c = cE(s) = E^p(s) =E^{p-1}(E(s)) = E^{p-1}(1)=0, 
\]
which implies $E^p=0$ and $E^{p-1}(u^{-1}) = u^{-1}c=0$.

By Lemma \ref{l additive slice}, we have
\[
R=\bigoplus_{i=0}^{p-1}R^E s^i.
\]
Write
$u^{-1}
=
\sum_{i=0}^{p-1}a_i s^i$ for some $a_0, ..., a_{p-1} \in R^E$. 
We have $0= E^{p-1}(u^{-1}) = (p-1)! a_{p-1}$, which implies $a_{p-1}=0$. Set
\[
t
:=
\sum_{i=0}^{p-2}
\frac{a_i}{i+1}s^{i+1} \in R.
\]
Since $a_i\in R^E$, we have
\[
E(t)
=E\left(\sum_{i=0}^{p-2}
\frac{a_i}{i+1}s^{i+1} \right)
=
\sum_{i=0}^{p-2}a_i s^i
=\sum_{i=0}^{p-1}a_i s^i
=
u^{-1}.
\]
Therefore, we obtain $D(t) = uE(t) = 1$. 
\end{proof}

\begin{dfn}\label{d nowhere vanishing}
Let $Y$ be a separated scheme of finite type over $k$. 
Let $D: \MO_Y \to \cL$ be a $k$-derivation, 
where $\cL$ is an invertible sheaf. 
We say that $D$ is {\em nowhere vanishing} if 
there exist an open cover 
$Y = \bigcup_{i \in I} Y_i$ and elements $t_i \in \Gamma(Y_i, \MO_Y)$  
such that $D(t_i)$ generates the $\MO_{Y_i}$-module $\cL|_{Y_i}$ for every $i \in I$. 
\end{dfn}

Let $Y$ be a variety and let $\cL$ be an invertible sheaf on $Y$. 
In general, a derivation $D: \MO_Y \to \cL$ does not necessarily induce $D^2$ nor $D^p$. 
On the other hand, if $\cL \simeq f^*\cL_X$ for some invertible sheaf $\cL_X$ on $Y^D =:X$, then 
we can naturally define the composition $D^n: \MO_Y \to \cL^{\otimes n}$ for every integer $n>0$.

\begin{prop}\label{p derivation to torsor}
Let $k$ be a field of characteristic $p>0$. 
Let $Y$ be a separated scheme of finite type over $k$ and let 
\[
D: \MO_Y \to \cL
\]
be a $k$-derivation, where $\cL$ is an invertible sheaf on $Y$. 
Let $f: Y \to Y^D =:X$ be the quotient by $D$. 
Assume that 
\begin{enumerate}
\item $D^p=0$, 
\item $D$ is nowhere vanishing, and 
\item $\cL \simeq f^*\cL_X$ for some invertible sheaf $\cL_X$ on $X$.  
\end{enumerate}
Then $f: Y \to X$ is an $\alpha_{\cL_X^{-1}}$-torsor. 
\end{prop}


\begin{proof}
We can find an affine open cover $X = \bigcup_{i \in I} X_i$ 
which trivialises $\cL_X$ and satisfies 
the condition in Definition \ref{d nowhere vanishing} for the pullback open cover $Y = \bigcup_{i \in I} Y_i$ and $t_i \in \MO_Y(Y_i)$. 
Hence we have 
$\MO_{X_i} \xrightarrow{\simeq, \varphi_i} 
 \MO_{X_i}e_i = \cL_X|_{X_i}$, and 
$\MO_{Y_i} \xrightarrow{\simeq, \psi_i} \MO_{Y_i}e_{i, Y} =\cL|_{Y_i}$   for 
some $e_i \in \cL_X(X_i)$ and $e_{i,Y}:=f^*e_i$. 
We have 
\[
D : \MO_{Y_i}
\xrightarrow{\partial_i} 
\MO_{Y_i} \xrightarrow{\psi_i, \simeq} \MO_{Y_i} e_{i, Y} = \cL|_{Y_i}, 
\]
where $\partial_i$ is a $k$-derivation. 
By $D^p=0$, we have $\partial_i^p=0$. 

The assumption on $t_i$ implies
that $\partial_i(t_i) \in \MO_{Y_i}(Y_i)^{\times}$. 
Therefore, after replacing $t_i$, 
Lemma \ref{l unit slice} allows us to assume that
$\partial_i(t_i)=1$ for every $i \in I$.  
Lemma \ref{l additive slice} now gives
\[
f_*\MO_{Y_i}
=
\bigoplus_{r=0}^{p-1}
\MO_{X_i}t_i^r \simeq  \MO_{X_i}[T_i]/(T_i^p - t_i^p). 
\]
Hence $f: Y \to X$ is locally an $\alpha_p$-torsor.

On the overlap $X_i \cap X_j$, we have 
$e_j=g_{ij}e_i$ in $\cL_X(X_i \cap X_j)$ for some $g_{ij}\in \MO_X(X_i \cap X_j)^\times$. Then it is easy to see that $\partial_j=g_{ij}^{-1}\partial_i.$ 
As $\partial_j(t_j)=1$, we get 
$\partial_i(t_j)=g_{ij}\partial_j(t_j) = g_{ij}$. 
This, together with $\partial_i(g_{ij})=0$ and $\partial_i(t_i)=1$, 
implies that 
\[
\partial_i(t_j-g_{ij}t_i) = \partial_i(t_j) - g_{ij} \partial_i(t_i) = g_{ij} - g_{ij} \cdot 1 =0.
\]
Therefore, $t_j=g_{ij}t_i+c_{ij}$ for some $c_{ij}\in\MO_X(X_i \cap X_j)$.
These are precisely the affine transition
functions of an $\alpha_{\cL_X^{-1}}$-torsor \cite[Page 38]{CDL25}, and hence 
$f: Y \to X$ is an $\alpha_{\cL_X^{-1}}$-torsor.
\end{proof}

\bibliographystyle{skalpha}
\bibliography{reference.bib}

\end{document}